\documentclass[a4paper, 12pt]{amsart}

\usepackage[font=small,labelfont=bf]{caption}
 \usepackage{amsmath, amsthm, amssymb, amsfonts}
\usepackage[utf8]{inputenc}
\usepackage[normalem]{ulem}
\usepackage[colorlinks = true,
            linkcolor = black,
            urlcolor  = blue,
            citecolor = red,
            anchorcolor = red]{hyperref}
\usepackage[nameinlink]{cleveref}
\usepackage{url}
\usepackage{comment}
\usepackage{mathrsfs}  
\usepackage{longtable} 
\usepackage{enumitem}
\usepackage{relsize}
\usepackage{bm}
\usepackage{xcolor}
\usepackage{graphicx}
\usepackage{mathabx}
\usepackage{mathdots}
\usepackage{mathtools}
\usepackage{color}
\usepackage{tikz-cd,tikz}
\usepackage{braket}
\usepackage{stmaryrd}
\usepackage{amssymb}
\usepackage{tipa}
\usepackage[margin=1.2 in]{geometry}
\usepackage{appendix}
\usepackage{physics}

\numberwithin{equation}{section}

\makeatletter
\newcommand{\leqnomode}{\tagsleft@true}
\newcommand{\reqnomode}{\tagsleft@false}
\makeatother

\usepackage{fancyhdr}
\usetikzlibrary{arrows, decorations.markings}
\usetikzlibrary{positioning}
\usetikzlibrary{matrix}
\usetikzlibrary{shapes}
\usetikzlibrary{arrows,decorations.markings, tikzmark}
\usepackage{stmaryrd}

\pgfarrowsdeclare{bad to}{bad to}
{
  \pgfarrowsleftextend{-2\pgflinewidth}
  \pgfarrowsrightextend{\pgflinewidth}
}
{
  \pgfsetlinewidth{0.8\pgflinewidth}
  \pgfsetdash{}{0pt}
  \pgfsetroundcap
  \pgfsetroundjoin
  \pgfpathmoveto{\pgfpoint{-3\pgflinewidth}{4\pgflinewidth}}
  \pgfpathcurveto
  {\pgfpoint{-2.75\pgflinewidth}{2.5\pgflinewidth}}
  {\pgfpoint{0pt}{0.25\pgflinewidth}}
  {\pgfpoint{0.75\pgflinewidth}{0pt}}
  \pgfpathcurveto
  {\pgfpoint{0pt}{-0.25\pgflinewidth}}
  {\pgfpoint{-2.75\pgflinewidth}{-2.5\pgflinewidth}}
  {\pgfpoint{-3\pgflinewidth}{-4\pgflinewidth}}
  \pgfusepathqstroke
}
\theoremstyle{definition}
\newtheorem{definition}{Definition}[section]

\newtheorem{theorem}[definition]{Theorem}
\newtheorem{example}[definition]{Example}
\newtheorem{proposition}[definition]{Proposition}

\newtheorem{corollary}[definition]{Corollary}

\newtheorem{lemma}[definition]{Lemma}

\newtheorem{question}{Question}
\newtheorem{remark}[definition]{Remark}

\newcommand{\BC}{{\mathbb{C}}}

\newcommand{\BG}{{\mathbb{G}}}

\newcommand{\BN}{{\mathbb{N}}}

\newcommand{\BZ}{{\mathbb{Z}}}

\newcommand\useleqno{\renewcommand\@eqnnum{\hb@xt@.01\p@{}%
                      \rlap{\normalfont\normalcolor
                        \hskip -\displaywidth(\theequation)}}}

\newcommand{\CE}{{\mathcal E}}

\newcommand{\CM}{{\mathcal M}}

\newcommand{\CT}{{\mathcal T}}

\DeclareMathOperator{\id}{id}

\DeclareMathOperator*{\res}{Res\,}

\newcommand{\Ext}{\mathcal{E}\text{xt}}

\renewcommand{\rank}{\mathop{\rm rank}\nolimits}

\newcommand{\End}{{\rm End}}

\newcommand{\GL}{{\rm GL}}

\newcommand{\pt}{{\mathsf{pt}}}

\newcommand{\congpf}{\xymatrix@1@=15pt{\ar[r]^-\sim&}}

\newcommand{\identity}{\mathrm{Id}}

\newcommand{\ch}{\mathrm{ch}}

\newcommand{\Perf}{\mathrm{Perf}}

\newcommand{\Hom}{\mathrm{Hom}}

\usepackage{tasks}

\usetikzlibrary{decorations.pathmorphing}
\AtBeginDocument{%
	\def\MR#1{}
}

\renewcommand{\Ext}{\textup{Ext}}

\newcommand{\normord}[1]{:\mathrel{#1}:}
\newcommand{\Yop}[1]{#1^{\text{op}}}
\newcommand{\modCQ}[1]{\mathrm{mod}\text{-}\mathbb{C}#1}

\usepackage{tikz}
\usepackage{tikz-cd}
\usetikzlibrary{decorations.pathmorphing}
\AtBeginDocument{%
	\def\MR#1{}
}

\begin{document}

\baselineskip=16pt
\parskip=5pt

\title[Field algebras associated with non-symmetric lattices]{Field algebras associated with non-symmetric lattices}
\author{Wacyl Meddour \\ Massachusetts Institute of Technology}

\email{meddour@mit.edu}

\date{11/21/2026}
\maketitle

\begin{abstract}
Field algebras, introduced by Bakalov and Kac, generalize vertex algebras by
replacing locality with a weaker associativity-type axiom. We develop the structure theory of field
algebras with a view toward geometric examples in which a natural bilinear form
fails to be symmetric. We prove a reconstruction theorem for field
algebras, producing a field algebra structure from a generating family of fields. We associate to every finite-rank lattice $L$ equipped with an
arbitrary bilinear form a canonical field algebra $V_L$, generalizing the Heisenberg and lattice vertex algebras. We then study potential definitions of commutants in this setting and analyze when they are realized as subfield algebras. Finally, we realize these constructions geometrically by looking at the homology of the moduli stack of complexes of quiver representations, extending previous results from Dominic Joyce.
\end{abstract}

\tableofcontents

\newpage
\section{Introduction} \label{sec: intro}

\subsection{Overview}

Vertex algebras \cite{Borcherds1986} were introduced as a generalization of commutative associative algebras and are most notable for their key role in Borcherds' proof of the monstrous moonshine conjecture, see \cite{Borcherds1992} \cite{Tatitscheff2019}. In essence, a vertex algebra consists of a vector space along with an integer-indexed collection of (not necessarily associative) distributive products. One should view the nonnegative-indexed products as generalizations of Lie brackets, while negative-indexed products correspond to a unital commutative algebra component. 

More formally, a vertex algebra consists of a vector space $V$ over a field $k$ along with a linear map $Y : V \rightarrow \End_k(V)\llbracket z,z^{-1}\rrbracket$ satisfying various conditions we describe in \Cref{sec: preliminaries}. Such a correspondence $Y$ associates to each vector $a \in V$ a series of endomorphisms $Y(a,z) =  \sum_{n \in \mathbb{Z}} a_{(n)}z^{-n-1}$, and such objects are often called fields, making $Y$ a state-field correspondence \cite{BakalovKac2003}. 

The main condition imposed on $Y$ is locality, which is effectively an analogue for commutativity. Bakalov and Kac \cite{BakalovKac2003} introduce the notion of \textit{field algebras} which replaces locality with a weaker condition analogous to associativity. The purpose of this paper is to build on the existing theory of field algebras motivated by geometric examples.

Vertex algebras already arise naturally from geometric examples. For example, one can associate a vertex (super)algebra to any finite-rank lattice $L$ equipped with a symmetric bilinear form \cite{FrenkelBenZvi2004}. Other examples include Joyce vertex algebras \cite{Joyce2019} which are vertex algebras associated to the homology of some $\mathbb{C}$-linear moduli stacks, such as representation spaces of quivers or sheaves on varieties.

In this paper, we aim at generalizing some of these techniques to the non-symmetric case. We begin by developing an analogue for the reconstruction theorem for field algebras. We then construct field algebras associated to arbitrary finite-rank lattices $L$ with non-symmetric pairings (free abelian groups equipped with a bilinear pairing which is not required to be symmetric). We then investigate commutants in field algebras, and conclude our discussion by extending Joyce's results \cite{Joyce2019} on the homology of the moduli stack of quiver representations to the non-symmetric case. 

\subsection{Main results}
This subsection summarizes the main results of our paper. 

\textbf{Theorem A}. (Reconstruction Theorem for field algebras)     Let $V$ be a vector space, $\ket{0} \in V$ and $T \in \End(V)$. Suppose we have a collection of fields $\{a^\alpha(z),b^\alpha(z)\}_{\alpha \in A}$ for some index set $A$. Suppose the following conditions hold: \begin{enumerate}
    \item For all $\alpha \in A$, $a^{\alpha}(z)\ket{0} \in a^\alpha + zV\llbracket z \rrbracket$ and $b^\alpha(z)\ket{0} \in a^\alpha + zV\llbracket z \rrbracket$ where $\{a^\alpha\}_{\alpha \in A} \subseteq V$ is a collection of vectors in $V$.
    \item $T\ket{0} = 0$, and $[T,a^\alpha(z)]= \partial_za^\alpha(z)$ and $[T,b^\alpha(z)] = \partial_z b^\alpha(z)$ for all $\alpha \in A$.
    \item There exists a collection $S$ of finite ordered tuples of elements of $A$ such that $V$ is spanned by  vectors of the form $a^{\alpha_1}_{(j_1)}a^{\alpha_2}_{(j_2)}\ldots a^{\alpha_m}_{(j_m)}\ket{0}$ for $j_i < 0$ and $(\alpha_1,\alpha_2,\cdots,\alpha_m) \in S$. 
    \item All pairs $(a^\alpha(z),b^\beta(z))$ are local for $\alpha,\beta \in A$. 
    \item For each tuple $(\alpha_1,\cdots,\alpha_m)\in S$, the pairs $(a^{\alpha_i}(z),a^{\alpha_j}(z))$ and $(b^{\alpha_i}(z),b^{\alpha_j}(z))$ are weakly local for any $1 \le i < j \le m$.
\end{enumerate} Then, the following assignment $$\begin{aligned} 
Y(a^{\alpha_1}_{(j_1)}a^{\alpha_2}_{(j_2)}\cdots a^{\alpha_m}_{(j_m)}\ket{0}, z) &= \dfrac{1}{(-j_1-1)!}\cdot \cdots \cdot \dfrac{1}{(-j_m-1)!} \normord{\partial_z^{-j_1-1}a^{\alpha_1}(z)\cdots\partial_z^{-j_m-1}a^{\alpha_m}(z)}
\end{aligned}$$  defines a field algebra structure on $V$, where $\Yop{Y}$ is defined exactly the same way as above with the fields $b^\alpha(z)$ instead of $a^\alpha(z)$. Furthermore, this is the unique field algebra structure on $V$ for which $Y(a^\alpha,z) = a^\alpha(z)$ with translation map $T$ and vacuum vector $\ket{0}$. In such a case, we say that $(V,\ket0, T,Y)$ is the \textit{field algebra generated by the collections} $\{a^\alpha(z)\}_{\alpha \in A}$ \textit{and} $\{b^\alpha(z)\}_{\alpha \in A}$.  \\

If $L$ is a finite-rank lattice with a symmetric pairing, one can construct a field algebra on: \[
V_L := \bigoplus_{\lambda \in L} \pi_\lambda
\] where each $\pi_\lambda$ is a copy of the Fock space, see \cite{FrenkelBenZvi2004}, chapter 5. This construction depends on some choice of elements $c_{\lambda,\mu} \in k$ for each $\lambda,\mu \in L$, usually taken in $\{-1,1\}$. If $L$'s pairing is not symmetric, a similar setup happens, in which we also have choices $c_{\lambda,\mu}^\text{op}$ for operators of $Y^\text{op}$. 

\textbf{Theorem B}. (Lattice field algebras). Let $L$ be a finite rank lattice with a bilinear form $(\cdot,\cdot) : L \times L \rightarrow \mathbb Z$. By setting $c_{\lambda,\nu} = 1$ and $\Yop{c}_{\lambda,\nu}=(-1)^{(\nu,\lambda)}$, the assignments \eqref{eq: exp-defn} and \eqref{eq: exp-defn-op} along with \eqref{eq: heisenberg-defn} extend to a field algebra structure on $V_L$ in the sense of \Cref{thm : reconstruction}. We call $(V_L,\ket 0, T,Y)$ the \textit{canonical field algebra structure associated to $L$}.\\

If $V$ is a vertex algebra and $v \in V$, we can define the commutant of $v$ as the set of elements $w \in V$ for which $v_{(n)}w=0$ for all $n \ge 0$, yielding a sub-vertex algebra of $V$ by OPE (see 2.2 in \cite{Kac1998}). However, this may fail in the general case of field algebras (see \Cref{ex: comm-not-subfield}). It remains true in the context of lattice field algebras for exponential elements:

\textbf{Theorem C}. (Commutant of a lattice exponential is a subfield algebra).     Let $L$ be a finite rank lattice with bilinear form $(\cdot,\cdot) : L \times L \rightarrow \BZ$, and let $(V_L,\ket 0, T,Y)$ be the (canonical) field algebra structure associated to $L$. Then for each $\alpha \in L$, the subspace \[
C(\alpha) := \bigcap_{n \in \BZ_{\ge 0}} \ker(\ket \alpha_{(n)}) \subseteq V_L
\] is a subfield algebra of $V_L$. \\

In \cite{Joyce2019}, Joyce associates a vertex algebra structure to the homology of the moduli stacks of objects in a triangulated category $\CT$ along with extra data of a perfect complex satisfying some compatibility conditions. The motivating picture for this project was the case where $\CT = D^b\modCQ{Q}$ for a quiver $Q$, in which case we get a lattice vertex algebra. The chosen perfect complex in that case is symmetrized. \\

We show that, by dropping the symmetrization condition, we recover the data of a lattice field algebra on the homology of the moduli stack of complexes of quiver representations. Furthermore, this naturally realizes the homology of the moduli stack of representations of $Q$ as an intersection of commutants if $Q$ is acyclic.

\textbf{Theorem D}. For any finite quiver $Q$, let $\overline{\CM}$ be the higher $\BC$-stack of complexes of representations for $Q$. There is an isomorphism of field algebras: \[
    \xi : H_*(\overline{\CM};\BC) \xrightarrow{\cong } V_{L_Q}
    \] realizing $H_*(\overline{\CM};\BC)$ as a lattice field algebra.

Furthermore, if $Q$ is an acyclic quiver and $\CM$ is the moduli stack of quiver representations, the map $\xi$ above identifies $H_*(\CM;\BC) \subset H_*(\overline{\CM};\BC)$ with an intersection of commutants of lattice exponentials, see \Cref{cor: moduli-stack-of-reps-is-commutant}. This is a continuation of previous work on the subject in the symmetric (vertex algebra) case, see \cite{Joyce2019}, \cite{Latyntsev2021}.

\subsection{Acknowledgments} The author thanks Miguel Moreira for his supervision, guidance and for many helpful conversations throughout this project, and for his useful comments. The author would also like to thank both Bojko Bakalov and Victor G. Kac for communications regarding the reconstruction theorem. This work was carried out as part of the UROP program, and the author thanks the UROP department for offering this opportunity.

\subsection{AI Disclosure} In the preparation of this paper, AI models were used to format the bibliography as well as to find references. The author takes full responsibility for the results presented in this paper.

\subsection{Remarks} In this paper, we assume that the underlying field $k$ of any field or vertex algebra discussed has characteristic $0$. We will also drop the mention of $k$ from now on.

\newpage 

\section{Preliminaries}  \label{sec: preliminaries}

We begin by recalling classical definitions of vertex and field algebras as well as some well-known properties of these objects. We will follow conventions in \cite{BakalovKac2003}. 

\subsection{State-field correspondences} 
\begin{definition}
    Let $V$ be a vector space. A \textit{field} is a series $f(z) \in \End(V)\llbracket z,z^{-1}\rrbracket$ written as: \[
    f(z) = \sum_{n \in \mathbb{Z}}f_{(n)}z^{-n-1}
    \] such that for each $b \in V$ there exists a threshold $N_b \in \BZ$ such that for all $n \ge N_b$ we have $f_{(n)} b = 0$.
\end{definition} 
Equivalently, a field $f$ satisfies $f(z)b \in V(\!(z)\!)$ for each $b \in V$, where $V(\!(z)\!)$ is the space of $V$-valued Laurent series.
\begin{definition}
A \textit{state-field correspondence} $(V,\ket{0}, T,Y)$ consists of the following data: \begin{itemize}
    \item A vector space $V$. 
    \item A vacuum vector $\ket{0} \in V$.
    \item A translation linear map $T \in \End(V)$
    \item A linear map $Y : V \rightarrow \End(V)\llbracket z,z^{-1}\rrbracket$ associating to each $a \in V$ a field $Y(a,z) = \sum_{n \in \mathbb{Z}} a_{(n)}z^{-n-1}$ satisfying:
    \begin{itemize}
        \item $Y(\ket{0}, z) = \identity_V$. 
        \item For all $a \in V$ we have $Y(a,z) \cdot \ket{0} \in V\llbracket z\rrbracket$ with $Y(a,z) \cdot \ket{0} |_{z=0} = a$, or equivalently $a_{(n)}\ket{0} = 0$ for $n \ge 0$ and $a_{(-1)}\ket{0} = a$.
        \item For all $a \in V$ we have $[T,Y(a,z)] = Y(Ta,z)= \partial_zY(a,z)$ or equivalently $[T,a_{(n)}] = -na_{(n-1)}$. 
    \end{itemize}
    Notice in particular this implies $T \ket{0} = 0$. 
\end{itemize}
\end{definition} 
The data of a state-field correspondence is enough to encapsulate integer-indexed products. They, however, may not be well-behaved in general. Vertex algebras require locality as an extra condition akin to commutativity. Field algebras require a weaker condition than locality akin to associativity. 
\begin{definition} \cite{BakalovKac2003}
    For every integer $n$, define the  $n$-\textit{th product} of two fields $a(z),b(z)$ as: \[
        a(z)_{(n)}b(z) = \res_x\big(a(x)b(z)i_{x,z}(x-z)^n-b(z)a(x)i_{z,x}(x-z)^n\big)
    \] where $i_{x,z}$ means expanding on $|x| > |z|$. Equivalently: \begin{equation}
    \label{eq: nth-prod}
        a(z)_{(n)}b(z) = \begin{cases}
            \res_x (x-z)^n[a(x),b(z)] & \text{if } n \ge 0 \\
            \normord{(\partial_z^{-n-1}a(z))b(z)}/(-n-1)! & \text{otherwise}
        \end{cases}
    \end{equation}
\end{definition}

In particular if $n \ge 0$ one has as coefficients: \begin{equation}
    \label{eq: nth-prod-coef}
    a(z)_{(n)}b(z) = \sum_{m \in \mathbb{Z}} \sum_{k=0}^n \binom{n}{k} (-1)^{n-k} [a_{(k)},b_{(m-k)}]z^{n-m-1}
\end{equation}
and for $n=-1$ we have: \[
a(z)_{(-1)}b(z) = \normord{a(z)b(z)} = a(z)_+b(z)+b(z)a(z)_-
\] where $a(z)_+ = \sum_{n < 0} a_{(n)}z^{-n-1}$ and $a(z)_- = a(z)-a(z)_+$. 

To each state-field correspondence $Y$, one can associate a dual correspondence $Y^{\text{op}}$. This will allow us to better formulate the necessary axioms for vertex and field algebras respectively later on. 

\begin{definition}
    Let $(V,\ket{0},T,Y)$ be a state-field correspondence. Define: \begin{equation}
     \label{eq: opp}
    \Yop{Y}(a,z)b = e^{zT}Y(b,-z)a       
    \end{equation} for all $a,b \in V$. In coefficients, this is equivalent to: \begin{equation}
        \label{eq: opp-coef}
        \sum_{n \in \BZ} \Yop{a_{(n)}}b\cdot z^{-n-1} = \sum_{n \in \mathbb{Z}} \sum_{k \ge 0} \dfrac{(-1)^{n+k+1}}{k!} T^k b_{(n+k)}a \cdot  z^{-n-1}
    \end{equation} where $\Yop{a}_{(n)}b$ is the $-n-1$'st coefficient of $\Yop{Y}(a,z)b$. Notice the right hand side is well-defined since $b_{(n)}a = 0$ for $n \gg 0$.
\end{definition} 

\begin{proposition}
    \label{prop: s2fprop}
    For any state-field correspondence, one has the following for all $a,b \in V$: \begin{enumerate}
        \item $Y(a,z)\ket{0} = e^{zT}a$
        \item $e^{wT}Y(a,z)e^{-wT} = Y(e^{wT}a,z) = i_{z,w}Y(a,z+w)$ 
        \item $\big (Y(a,z)_{(n)}Y(b,z) \big )\ket{0} = Y(a_{(n)}b,z)\ket{0}$ for all $n \in \mathbb{Z}$
        \item $(V,\ket{0},T,\Yop{Y})$ is a state-field correspondence and $(Y^{\text{op}})^{\text{op}} = Y$
    \end{enumerate}
\end{proposition} 
\begin{proof}
    Propositions 1.7 and 1.8 in \cite{BakalovKac2003}. 
\end{proof}
\subsection{Vertex and field algebras} We now discuss vertex and field algebras. As explained above, one must attach additional conditions to state-field correspondences to obtain these objects. 

\begin{definition} Consider two fields $a(z),b(z)$ on a vector space $V$. 
    \begin{itemize}
        \item We say that the pair $(a(z),b(z))$ is \textit{local} if for $N \gg 0$ one has:\begin{equation}
            \label{eq: locality}
            (z-w)^N[a(z),b(w)]=  0
        \end{equation} In coefficients, this is equivalent to: \begin{equation}
            \label{eq: locality-coef}
            \sum_{k=0}^N \binom{N}{k}(-1)^{N-k}[a_{(n+k)},b_{(m-k)}] = 0
        \end{equation} for all $n,m \in \mathbb{Z}$. ~\\ 
        \item We say that the pair $(a(z),b(z))$ is \textit{weakly local} if $a(z)_{(N)}b(z) = 0$ for $N \gg 0$. In coefficients, this is equivalent to: 
        \begin{equation}
            \label{eq: weak-locality-coef}
            \sum_{k=0}^N \binom{N}{k}(-1)^{N-k}[a_{(k)},b_{(m-k)}] = 0
        \end{equation} for all $m \in \mathbb{Z}$
    \end{itemize}
\end{definition} 

\begin{remark}
    Locality of $(a(z),b(z))$ implies weak locality of $(a(z),b(z))$.
\end{remark}
Notice that weak locality is not symmetric, but locality is. We are now ready to define vertex and field algebras. 
\begin{definition} Consider a state-field correspondence $(V,\ket{0},T,Y)$ \begin{itemize}
    \item We say that $Y$ is a \textit{vertex algebra} if all pairs $(Y(a,z),Y(b,z))$ are local for $a,b\in V$. This property is called \textit{locality}.
    \item We say that $Y$ is a \textit{field algebra} if all pairs $(Y(a,z),\Yop{Y}(b,z))$ are local for $a,b\in V$.
\end{itemize}
\end{definition}

A state-field correspondence $(V,\ket 0, T,Y)$ is \textit{associative} if for all $a,b,c \in V$ the following relation holds for all $N \gg 0$:
    \begin{equation}
    \label{eq: associativity}
    (z-w)^NY(Y(a,z)b,-w)c = (z-w)^N i_{z,w} Y(a,z-w)Y(b,-w)c
\end{equation} where \[
Y(Y(a,z)b,-w) = \sum_{n \in \mathbb{Z}}Y(a_{(n)}b,-w)z^{-n-1}
\] This relation is called \textit{associativity}. 

In coefficients, associativity is equivalent to: \begin{equation}
    \label{eq: associativity-coef}
    \sum_{k=0}^N \binom{N}{k} (a_{(n+k)}b)_{(m-k)}c = \sum_{k=0}^N \sum_{j \ge 0} \binom{N}{k} \binom{-n-k+j-1}{j} a_{(n+k-j)} b_{(m+j-k)} c
\end{equation} for all $n,m \in \mathbb{Z}$.

In \cite{BakalovKac2003}, field algebras are defined as associative state-field correspondences. The two definitions are equivalent:

\begin{proposition}
    \label{prop: characterization}
    Fix a vector space $V$. \begin{enumerate}
        \item A state-field correspondence $Y$ on $V$ is associative if and only if there exists a state-field correspondence $X$ for which all pairs $(Y(a,z),X(b,z))$ are local for $a,b \in V$, in which case  $X=\Yop{Y}$. 
        \item Any vertex algebra is a field algebra.
        \item A field algebra $Y$ on $V$ is a vertex algebra if and only if $Y=\Yop{Y}$. 
    \end{enumerate}
\end{proposition} 
\begin{proof}
    (1) is Theorem 3.7 in \cite{BakalovKac2003} along with Proposition 2.10. (2) follows from (1) and Proposition 2.10. (3) is Theorem 6.3. 
\end{proof} 
In fact, it turns out that any vertex algebra is a \textit{strong field algebra}, i.e a state-field correspondence satisfying the $n$-\textit{th product axiom} \cite{BakalovKac2003}: \[
Y(a,z)_{(n)}Y(b,z)=Y(a_{(n)}b,z)
\] This is a consequence of OPE expansion \cite{FrenkelBenZvi2004}. There is also an equivalent condition for strong field algebras: They are field algebras in which all pairs $(Y(a,z),Y(b,z))$ are weakly local. 

\begin{remark}
    In \cite{Kac1998}, field algebras are defined as what we have called strong field algebras. This definition was corrected in \cite{BakalovKac2003} after the realization that requiring weak locality of all pairs was too strong.  
\end{remark}

Although we will require weak locality in our reconstruction theorem for field algebras, we will not need it to hold for all such pairs, so we will not discuss the concept of strong field algebras.

\begin{proposition}
    \label{prop : ope-goddard}
    Let $(V,\ket{0},T,Y)$ be a state-field correspondence. Let $a(z),b(z)$ be two fields on $V$. \begin{enumerate}
        \item (OPE) The pair $(a(z),b(z))$ is local if and only if there exists $K \ge 0$ for which: \[
        [a(z),b(w)] = \sum_{j=0}^K \big ( a(w)_{(j)}b(w)\big ) \cdot \dfrac{\partial_w^j  \delta(z-w)}{j!}
        \] where $\delta(z-w) = \sum_{n \in \mathbb{Z}} z^nw^{-n-1}$. In particular if $Y$ is a vertex algebra we get OPE expansion: $$[Y(a,z),Y(b,w)] = \sum_{j=0}^K \dfrac{Y(a_{(j)}b,w)}{j!} \cdot \partial_w^j \delta(z-w)$$
        \item (Goddard's Uniqueness) If $a(z)\ket{0} = b(z)\ket{0}$ and all pairs $(Y(v,z),a(z))$ and $(Y(v,z),b(z))$ are local for $v \in V$, then $a(z)=b(z)$. 
    \end{enumerate}
\end{proposition} \begin{proof}
    See 2.2 in \cite{Kac1998} and 2.4 in \cite{BakalovKac2003} respectively. 
\end{proof}
 Dong's Lemma allows us to understand how locality behaves under $n$-th products of fields. We shall use it in the proof of the reconstruction theorem for field algebras. This is Lemma 2.5 in \cite{BakalovKac2003}.

\begin{lemma} (Dong's Lemma)
    \label{lemma: dong}
    Let $a,b,c$ be fields over $V$:\begin{enumerate}
        \item If $(a,b)$ is weakly local and $(a,c),(b,c)$ are local then $(a_{(n)}b,c)$ is local for any $n \in \mathbb{Z}$. 
        \item If $(a,b)$ and $(a,c)$ are weakly local and $(b,c)$ is local then $(a_{(n)}b,c)$ is local for $n \ge 0$.
        \item If all pairs $(a,b)$, $(a,c)$, $(b,c)$ are weakly local then $(a_{(n)}b,c)$ is weakly local for any $n \ge 0$ and $(a,b_{(n)}c)$ is weakly local for all $n \in \mathbb{Z}$.
    \end{enumerate}
\end{lemma} \begin{proof}
    See Section 3.2 in \cite{Kac1998}. 
\end{proof}

We conclude this section by a characterization of weak locality, similar to Corollary 2.2 in \cite{Kac1998}.

\begin{lemma}
    \label{lemma: weak-loc-formulation}
    Consider a formal distribution over a vector space $U$: \[
    f(z,w) = \sum_{n,m \in \mathbb{Z}} A_{n,m} \cdot z^{n}w^{m} \in U\llbracket z,z^{-1},w,w^{-1}\rrbracket
    \] Let $N \ge 0$. Then $\text{Res}_z (z-w)^n f(z,w)=0$ for all $n \ge N$ if and only if we can find $\psi_j(w) \in U\llbracket w,w^{-1}\rrbracket$ with $j=0,\ldots,N-1$ such that: \[
    f(z,w)_{z_-} = \sum_{j=0}^{N-1}\psi_j(w)\partial_w^j\delta(z-w)_{z_-}
    \] where $f(z,w)_{z_-} = \displaystyle\sum_{n < 0} \sum_{m \in \mathbb Z} A_{n,m} z^nw^m $. 
\end{lemma}
\begin{proof}
    For all $n \ge 0$ one has $\text{Res}_z (z-w)^n f(z,w) = \text{Res}_z (z-w)^n f(z,w)_{z_-}$, and one can check that $\text{Res}_z(z-w)^n \partial_w^j \delta(z-w)_{z_-}=0$ for $n \ge j+1$. So the result follows by Proposition 2.2(c) in \cite{Kac1998}. 
\end{proof}

Weak locality of $(a(z),b(z))$ is then the condition of \Cref{lemma: weak-loc-formulation} with $f(z,w) = [a(z),b(w)]$. 

\newpage

\section{Reconstruction theorem for field algebras}  \label{sec: reconstruction}

\subsection{Reconstruction for field algebras}
Recall the reconstruction theorem for vertex algebras as described in Section 4.4 in \cite{FrenkelBenZvi2004}. Given a collection of local ``top'' fields generating our vector space $V$, one can uniquely reconstruct a vertex algebra structure on $V$. An analogous result with essentially the same proof exists for field algebras, provided we require the generating fields to be weakly local. 

As \Cref{prop: characterization} suggests, we will interpret field algebras as local pairs of state-field correspondences $Y,\Yop{Y}$, since locality is easier to grasp than associativity as a condition. Throughout this section, analogous results hold by flipping the roles of $Y$ and $\Yop{Y}$. 


\begin{theorem}
    \label{thm : reconstruction}
    Let $V$ be a vector space, $\ket{0} \in V$ and $T \in \End(V)$. Suppose we have a collection of fields $\{a^\alpha(z),b^\alpha(z)\}_{\alpha \in A}$ for some index set $A$. Suppose the following conditions hold: \begin{enumerate}
        \item For all $\alpha \in A$, $a^{\alpha}(z)\ket{0} \in a^\alpha + zV\llbracket z\rrbracket$ and $b^\alpha(z)\ket{0} \in a^\alpha + zV\llbracket z\rrbracket$ where $\{a^\alpha\}_{\alpha \in A} \subseteq V$ is a collection of vectors in $V$.
        \item $T\ket{0} = 0$, and $[T,a^\alpha(z)]= \partial_za^\alpha(z)$ and $[T,b^\alpha(z)] = \partial_z b^\alpha(z)$ for all $\alpha \in A$.
        \item There exists a collection $S$ of finite ordered tuples of elements of $A$ such that $V$ is spanned by  vectors of the form $a^{\alpha_1}_{(j_1)}a^{\alpha_2}_{(j_2)}\ldots a^{\alpha_m}_{(j_m)}\ket{0}$ for $j_i < 0$ and $(\alpha_1,\alpha_2,\cdots,\alpha_m) \in S$. 
        \item All pairs $(a^\alpha(z),b^\beta(z))$ are local for $\alpha,\beta \in A$. 
        \item For each tuple $(\alpha_1,\cdots,\alpha_m)\in S$, the pairs $(a^{\alpha_i}(z),a^{\alpha_j}(z))$ and $(b^{\alpha_i}(z),b^{\alpha_j}(z))$ are weakly local for any $1 \le i < j \le m$.
    \end{enumerate} Then, the following assignment $$\begin{aligned} 
Y(a^{\alpha_1}_{(j_1)}a^{\alpha_2}_{(j_2)}\ldots a^{\alpha_m}_{(j_m)}\ket{0}, z) &= \dfrac{1}{(-j_1-1)!}\cdot \cdots \cdot \dfrac{1}{(-j_m-1)!} \normord{\partial_z^{-j_1-1}a^{\alpha_1}(z)\cdots\partial_z^{-j_m-1}a^{\alpha_m}(z)} \\
&= a^{\alpha_1}(z)_{(j_1)}( a^{\alpha_2}(z)_{(j_2)} \cdots (a^{\alpha_m}(z)_{(j_m)}(\text{Id}_V) )\cdots)
\end{aligned}$$  defines a field algebra structure on $V$, where $\Yop{Y}$ is defined exactly the same way as above with the fields $b^\alpha(z)$ instead of $a^\alpha(z)$. Furthermore, this is the unique field algebra structure on $V$ for which $Y(a^\alpha,z) = a^\alpha(z)$ with translation map $T$ and vacuum vector $\ket{0}$. In such a case, we say that $(V,\ket0, T,Y)$ is the \textit{field algebra generated by the collections} $\{a^\alpha(z)\}_{\alpha \in A}$ \textit{and} $\{b^\alpha(z)\}_{\alpha \in A}$. 
\end{theorem}
\begin{remark}
    Notice that we only require weak locality in a certain order. This is not equivalent to requiring weak locality of all pairs of generators, since weak locality is not symmetric. A similar version of the above theorem was provided in \cite{Kac1999} (the erratum of \cite{Kac1998}) but concluded that the resulting structure is that of a strong field algebra (a field algebra in the sense of \cite{Kac1998}), which is incorrect. See \Cref{ex: heisenberg}.
\end{remark}
\begin{proof}
    Fix any basis of $V$ of elements $a^{\alpha_1}_{(j_1)}\cdots a^{\alpha_m}_{(j_m)}\ket 0$ with $j_i < 0$ and $(\alpha_1,\cdots,\alpha_m) \in S$. Then as in Theorem 4.4.1 in \cite{FrenkelBenZvi2004}, we recover two state-field correspondences $Y$ (with $a^\alpha$'s) and $X$ (with $b^\alpha$'s) with the above assignment. 

    We begin by showing $Y$ and $X$ are local. By \Cref{lemma: dong} (3), we inductively see that $(b^{\alpha_1}(z), b^{\alpha_2}(z)_{(j_2)}(b^{\alpha_3}(z)_{(j_3)}(\ldots)))$ is weakly local for each generator $a^{\alpha_1}_{(j_1)}a^{\alpha_2}_{(j_2)}\ldots a^{\alpha_m}_{(j_m)}\ket{0}$ in our basis. Notice furthermore that locality (resp. weak locality) of $(a(z),b(z))$ remains when differentiating either side by \Cref{prop : ope-goddard} and \Cref{lemma: weak-loc-formulation}. As a consequence, any field $k(z)$ local with the generating fields $b^{\alpha}(z)$ with $\alpha \in A$ will be local with $X$ by \Cref{lemma: dong} (1), which implies $X$ and $Y$ are local so that \Cref{prop: characterization} gives $X=\Yop{Y}$ and hence a field algebra structure on $V$. 

    If $Y'$ is another field algebra on $V$ agreeing with $Y$ on top vectors $a^\alpha$ with $\alpha \in A$, the above shows that $X'=\Yop{(Y')}$ is local with $Y$ forcing $X'=\Yop{Y}=X$ and hence $Y=\Yop{X}=Y'$ again by \Cref{prop: characterization}. 
\end{proof}  
\begin{corollary}
    \label{cor: locality-check}
    If $(V,\ket0,T,Y)$ is a field algebra generated by $\{a^\alpha(z)\}_{\alpha \in A}$ and $\{b^\alpha(z)\}_{\alpha \in A}$, then a field $k(z) \in \End(V)\llbracket z,z^{-1}\rrbracket$ is local with $Y$ if and only if it is local with $\{a^\alpha(z)\}_{\alpha \in A}$. 
\end{corollary}

Weak locality of the generating fields does not necessarily carry down to the whole space. The Heisenberg field algebra will provide a counter-example. If, however, we do have weak locality, the $n$-th product relation holds:
\begin{proposition}
    \label{prop: topfields_nprod}
    Suppose $(V,\ket{0},T,Y)$ is a field algebra. If $(Y(a,z),Y(b,z))$ is weakly local for $a,b \in V$, then: \[
    Y(a,z)_{(n)}Y(b,z) = Y(a_{(n)}b,z)
    \] for all $n \in \mathbb{Z}$. 
\end{proposition}
\begin{proof}
    By \Cref{prop: s2fprop}, we see that both sides of the equation agree on the vacuum. By \Cref{lemma: dong} (1), both sides are local with $\Yop{Y}$, so we are done by \Cref{prop : ope-goddard} (2). 
\end{proof}

As a consequence, the $n$-th product axiom does hold for generating fields in the context of the reconstruction theorem. This result can be interpreted as motivation for requiring weak locality as an assumption in \Cref{thm : reconstruction}, as if $(\alpha_1,\alpha_2) \in S$, the expression $Y(a^{\alpha_1},z)_{(n)}Y(a^{\alpha_2},z)$ is precisely how we define $Y(a^{\alpha_1}_{(n)}a^{\alpha_2},z)$. 

Recall that weak-locality is not symmetric. In \Cref{thm : reconstruction}, we assume that only certain pairs of generating fields are weakly local. In the Heisenberg field algebra example we will present, all generating fields will be pairwise weakly local. That may not be true for field algebras built from \Cref{thm : reconstruction} in general, and the lattice field algebra will provide a counterexample.

\begin{proposition}
    \label{prop: rec-dual-prop}
    Suppose that $(V,\ket0, T,Y)$ is a field algebra generated by fields $\{a^\alpha(z)\},\{b^\alpha(z)\}$, $\alpha \in A$. Let $\alpha \in A$ such that $(a^\alpha(z),a^\beta(z))$ is weakly local for all $\beta \in A$. Then, for all $n \ge 0$ and $\beta \in A$: \[
    a^{\alpha}(z)_{(n)}b^\beta(z) = Y(a^\alpha,z)_{(n)}\Yop{Y}(a^\beta,z) = \Yop{Y}(a^\alpha_{(n)}a^\beta,z)
    \]
\end{proposition}
\begin{proof}
    Set $B_1(z) = a^\alpha(z)_{(n)}b^\beta(z)$ and $B_2(z) = \Yop{Y}(a^\alpha_{(n)}a^\beta,z)$. By \Cref{lemma: dong} (2), we see that for any $\gamma \in A$, $(B_1(z),a^\gamma(z))$ is local. Hence, as a consequence of the proof of the reconstruction theorem, we see that $B_1(z)$ is local to $Y$ by \Cref{cor: locality-check}. 

    The last thing to check before invoking \Cref{prop : ope-goddard} (2) is that $B_1(z) \ket 0 = B_2(z) \ket 0$. Denote $b^\beta(z) = \displaystyle \sum_{n \in \mathbb{Z}} b^\beta_{(n)} z^{-n-1}$. Recall by \cref{prop: s2fprop} we have that: \[
    B_2(z)\ket 0 = e^{zT} a^{\alpha}_{(n)}a^\beta = \sum_{m \ge 0} \dfrac{T^m}{m!} a^{\alpha}_{(n)}a^\beta z^m
    \]
    On the other hand, by \eqref{eq: nth-prod} we have that: \[
    \begin{split}
    B_1(z)\ket 0 &= \res_x (x-z)^n [a^\alpha(x),b^\beta(z)] \ket 0 \\
    &= \res_x \sum_{t=0}^n \sum_{k,l \in \mathbb Z} \binom{n}{t}(-1)^{n-t}[a^{\alpha}_{(k)},b^{\beta}_{(l)}]\ket 0 x^{t-k-1}z^{n-t-l-1} \\ 
    &= \sum_{t=0}^n \sum_{l \in \mathbb Z}  \binom{n}{t} (-1)^{n-t} [a^{\alpha}_{(t)},b^{\beta}_{(l)}] \ket 0 z^{n-t-l-1} \\
    &= \sum_{l \in \mathbb Z} \sum_{t=0}^n \binom{n}{t} (-1)^{t} [a^{\alpha}_{(n-t)}, b^\beta_{(t-l-1)}]\ket 0 z^l 
    \end{split}
    \]
    Notice $a^\alpha_{(n-t)}\ket 0 = 0$ as $n-t \ge 0$. If $l < 0$ then $t-l-1 \ge 0$ so the term vanishes. Using the fact $a^\beta(z)\ket 0 = e^{zT}a^\beta =  b^\beta(z)\ket 0$, we see that this identity boils down to: \[
    \dfrac{T^m}{m!} a^{\alpha}_{(n)}a^\beta = \sum_{t=0}^n \binom{n}{t} (-1)^{t} a^{\alpha}_{(n-t)} a^\beta_{(t-m-1)}\ket 0
    \] for all $m \ge 0$. This is true by induction on $m$ since $[T,a^\gamma_{(k)}] = -ka^\gamma_{(k-1)}$. 
\end{proof}
As a corollary of the above and \Cref{prop : ope-goddard}, we see that OPE (as usually done in vertex operator algebras) holds on ``source'' top fields:
\begin{corollary}
    \label{cor: ope-on-source-fields}
    Suppose that $(V,\ket0, T,Y)$ is a field algebra generated by fields $\{a^\alpha(z)\},\{b^\alpha(z)\}$, $\alpha \in A$. Let $\alpha \in A$ such that $(a^\alpha(z),a^\beta(z))$ is weakly local for all $\beta \in A$. Then, for all $\beta \in A$: \[
    [Y(a^\alpha,z),\Yop{Y}(a^\beta,w)] = \sum_{\substack{j \ge 0 \\ j \text{ finite}}} \dfrac{\Yop{Y}(a^\alpha_{(j)}a^\beta,w)}{j!} \partial_w^j \delta(z-w)
    \]
\end{corollary}

\subsection{Application to the Heisenberg case} In this section, we associate a field algebra analogue to the Heisenberg vertex algebra for any finite-rank lattice $L$. Symmetry of the pairing in the vertex algebra case is required due to the natural antisymmetry of the commutator. To fix this, we ``force'' our commutator to prefer one sign over the other. Since we have two possible choices of signs, this gives a natural definition of $\Yop{Y}$. We mimic the conventions of Section 5.4.1 in \cite{FrenkelBenZvi2004}. 

Begin with a finite-rank lattice $L$ with a (not necessarily symmetric) bilinear form $( \cdot ,\cdot) : L \times L \rightarrow \mathbb{Z}$, i.e $L$ is a finite-rank free abelian group with pairing $(\cdot,\cdot)$ not required to be symmetric. We complexify our lattice $L$ to get $\mathfrak{h} = L \otimes_{\mathbb{Z}} \mathbb{C}$. Consider the following central extension: \[
0 \rightarrow \mathbb{C} \textbf{1} \rightarrow \hat{\mathfrak{h}} \rightarrow \mathfrak{h}((t)) \rightarrow 0
\] with commuting relations defined by: \[[A \otimes f(t), B \otimes g(t)] = -(A,B)\text{Res}_t(f_+dg)\cdot \textbf{1} - (B,A) \text{Res}_t (f_-dg)\cdot \textbf{1} \] where if $f(t) = \sum_n a_n t^n$ then $f_+(t) = \sum_{n \ge 0} a_nt^n$ and $f_-(t) = f(t)-f_+(t)$. 

The above relations define a valid two-cocycle satisfying the Jacobi identity by centrality of $\textbf{1}$. Furthermore, we clearly have antisymmetry since: \[
\text{Res}_t(f_+dg) = -\text{Res}_t (g_-df)
\] The main relation we aim to obtain for our Heisenberg field algebra is \begin{equation}
    \label{eq: heisenberg-comm}
    [h_n,g_m] = n\delta_{m+n,0} \cdot ( (h,g) \textbf{1}_{n\ge 0} +(g,h)\textbf{1}_{n < 0})
\end{equation} where $\textbf{1}_{n \ge 0}$ is the indicator of the event that $n \ge 0$ (and likewise for $\textbf{1}_{n < 0}$). One can obtain this by taking the Weyl algebra $\Tilde{\mathcal{H}} = \Tilde{\mathcal{H}}_L $ as the completed enveloping algebra of $\hat{\mathfrak{h}}$ quotiented by $1 = \textbf{1}$. It has topological generators $h_n$, $h \in \mathfrak{h}$, $n \in \mathbb{Z}$ satisfying \eqref{eq: heisenberg-comm} where $h_n = ht^n$.  Finally, define the Fock space $\pi$ as usual by $\pi = \tilde{\mathcal{H}} \otimes_{\tilde{\mathcal{H}}_+} \mathbb{C}$, where $\tilde{\mathcal{H}}_+$ is the subalgebra of $\Tilde{H}$ generated by $h_n$ for $n \ge 0$, $h \in \mathfrak h$. 

More concretely, we think of the Fock space $\pi$ as polynomials in $h_{-n-1}$ with $h \in \mathfrak{h}$ and $n \ge 0$, where negative operators act by product and nonnegative ones act via \eqref{eq: heisenberg-comm} and $h_n \cdot 1 = 0$. We conflate $1 = \ket 0$.

The translation operator $T : \pi \rightarrow \pi$ is defined by the relations: \[
\begin{split}
    T \cdot \ket 0 &= 0 \\
    [T,h_n] &= -nh_{n-1} \quad \text{ for all } h\in \mathfrak{h} \text{ and } n< 0
\end{split}
\]
The second relation readily extends to $[T,h_n]=-nh_{n-1}$ for all $n \in \mathbb Z$. 

To each $h \in \mathfrak{h}$, associate a field to $h_{-1}\ket0$ defined by: \begin{equation}
    \label{eq: heisenberg-defn}
    h(z) = \sum_{n \in \mathbb{Z}} h_nz^{-n-1}
\end{equation}
For the reconstruction theorem to work, we will need the following to hold:
\begin{proposition}
    \label{prop : heisenberg-weakloc}
    The pair of fields $(h(z),g(z))$ is weakly local for all $h,g \in \mathfrak{h}$. 
\end{proposition}
\begin{proof}
    Follows from the following calculation for $N \ge 0$ and for all $m \in \mathbb{Z}$: \[
    \begin{split}
        \sum_{k=0}^N \binom{N}{k}(-1)^{N-k}[h_k,g_{m-k}] &= \textbf{1}_{m=0}\cdot  (h,g)\cdot \sum_{k=0}^N \binom{N}{k}(-1)^{N-k}\cdot k\\
        &=\textbf{1}_{m=0}\cdot  (h,g)\cdot \bigg ( \dfrac{d}{dx}(x-1)^N\bigg )|_{x=1} \\
        &= 0
    \end{split}
    \] for $N \ge 2$. 
\end{proof}

To show that our construction will yield a field algebra, we will invoke the reconstruction theorem. For this, we will need generating fields for the opposite state-field correspondence. Due to the construction of our Fock space $\pi$, it is enough to define an operator on $\pi$ by defining its action on $\ket{0}$ and how it commutes with each $h_i$ for $i < 0$.

\begin{definition} Fix $h \in \mathfrak{h}$ and $n \in \mathbb{Z}$. Define \textit{opposite operators} $\Yop{h_n} : \pi \rightarrow \pi$ for each $n \in \mathbb{Z}$ by the relations: \begin{equation}
\label{eq : heisenberg-op}
 \begin{split}
    \Yop{h_n} \ket 0 &= h_n \ket 0 \\
    [\Yop{h_n},g_m] &= n(g,h)\delta_{m+n,0}
\end{split}   
\end{equation}
for all $g \in \mathfrak{h}$ and $m \in \mathbb{Z}_{< 0}$. The second relation extends to all $m \in \mathbb Z$. We also define \textit{opposite fields}: \[\Yop{h}(z) = \sum_{n \in \mathbb{Z}}\Yop{h_n}z^{-n-1}\]
\end{definition}
Notice that $\Yop{h}_n=h_n$ for $n < 0$.

The opposite fields will generate $\Yop{Y}$ in our construction. As one would expect, the definitions could have been flipped, so the following commutation relation holds for opposite operators: 
\begin{lemma}
    \label{lemma: heisenberg-op-commute}
    For any $h,g \in \mathfrak{h}$ and $n,m \in \mathbb{Z}$ we have that: \[
    [\Yop{h_n}, \Yop{g_m}] = n\delta_{m+n,0} \cdot ((g,h) \textbf{1}_{n \ge 0} + (h,g) \textbf{1}_{n < 0} ) 
    \] In other words, $[\Yop{h_n}, \Yop{g_m}]  = [g_n,h_m]$
\end{lemma}
\begin{proof}
    Any operator on $\pi$ is defined by how it acts on $\ket 0$ and how it commutes with each $k_l$, $k \in \mathfrak{h}$, $l < 0$. By the Jacobi identity, for each $l \in \mathbb{Z}_{< 0}$: \[
    \begin{split}
        [[\Yop{h_n},\Yop{g_m}],k_l] &= [\Yop{h_n},[\Yop{g_m},k_l]]+[\Yop{g_m},[k_l,\Yop{h_n}]] \\
        &= [\Yop{h_n},m\delta_{m+l,0}(k,g)]+[\Yop{g_m},-n\delta_{n+l,0}(k,h)]\\
        &= 0
    \end{split}
    \]
    A calculation shows $[\Yop{h_n},\Yop{g_m}] \ket 0 =  n\delta_{m+n,0} \cdot ((g,h) \textbf{1}_{n \ge 0} + (h,g) \textbf{1}_{n < 0} )$, which concludes. 
\end{proof}
Opposite operators will play the role of $\Yop{Y}$ in our reconstruction theorem. Recall that we defined: \[
\Yop{h}(z) = \sum_{n \in \mathbb{Z}} \Yop{h_n}z^{-n-1}
\]
Notice that $\Yop{h}_{-1}\ket 0 = h_{-1}\ket 0$ and a similar argument to the above shows that $[T,\Yop{h}_n] = -n\Yop{h_{n-1}}$ for all $n \in \mathbb{Z}$. 
\begin{proposition}
    \label{prop: heisenberg-opp-locality}
    Let $h,g \in \mathfrak{h}$. Then the pair $(\Yop{h}(z),g(z))$ is local and the pair $(\Yop{h}(z),\Yop{g}(z))$ is weakly local.
\end{proposition} 
\begin{proof}
    The first part follows from the following calculation and \Cref{prop : ope-goddard}: \[
    \begin{split}
        [\Yop{h}(z),g(w)] &= \sum_{n,m \in \mathbb{Z}} [\Yop{h_n},g_m]z^{-n-1}w^{-m-1}\\
        &= -(g,h) \partial_z \sum_{n \in \mathbb Z} z^{n}w^{-n-1} \\
        &= -(g,h)\partial_z \delta(z-w)
    \end{split}
    \]
    The second part is similar to \Cref{prop : heisenberg-weakloc}. 
\end{proof} 
We thus have all components for the reconstruction theorem to apply:
\begin{corollary}
    The Fock space $\pi$ defined above has a natural structure of a field algebra $(\pi,\ket 0, T, Y)$ generated by $\{h(z)\}_{h \in \mathfrak{h}}$ and $\{\Yop{h}(z)\}_{h \in \mathfrak{h}}$. To each vector $a= h^{(1)}_{j_1}\ldots h^{(n)}_{j_n}\ket 0$ with $j_i < 0$, it assigns: \[
    Y(a,z) = h^{(1)}(z)_{(j_1)}(h^{(2)}(z)_{(j_2)}(\ldots h^{(n)}(z))\ldots )
    \] and similarly for $\Yop{Y}(a,z)$ with $\Yop{(h^{(i)})}(z)$ instead. 
\end{corollary}
\begin{example}
    If the bilinear form $(\cdot,\cdot) : L \times L \rightarrow \mathbb Z$ is symmetric, we automatically see that $h_n = \Yop{h_n}$ which implies $Y=\Yop{Y}$ so we recover the usual Heisenberg vertex algebra structure by \Cref{prop: characterization}. 
\end{example}

\begin{example}
\label{ex: heisenberg}
This construction is not a strong field algebra in general. Pick $L = \mathbb{Z}a \oplus \mathbb{Z}b$ with pairing given by \[Q=
\begin{pmatrix}
    ( a,a ) & ( a,b) \\
    ( b,a) & ( b,b)
\end{pmatrix} = \begin{pmatrix}
    0 & -1 \\
    1 & 0
\end{pmatrix}
\] Consider the Heisenberg field algebra for this lattice. Then: \[
Y(a_{(-1)}a_{(-1)}\ket 0,z)=\sum_{n \in \mathbb{Z}} \sum_{k+l = n-1} a_{(k)}a_{(l)} z^{-n-1}
\] since $( a,a) = 0$. If we call the inner coefficient $l_{(n)}$ as a well-defined operator on our vector space, then for all $N \ge 0$: \[
\begin{split}
    Y(a_{(-1)}a_{(-1)}\ket 0,z)_{(N)}Y(b_{(-1)}\ket 0,z) &= \res_x (x-z)^N [Y(a_{(-1)}a_{(-1)}\ket 0,x), Y(b_{(-1)}\ket 0, z)] \\
    &= \res_x \sum_{t=0}^N \sum_{n,m} \binom{N}{t}(-1)^{N-t} [l_{(n)},b_{(m)}]x^{t-n-1}z^{N-t-m-1} \\
    &= \sum_{t=0}^N \sum_m \binom{N}{t}(-1)^{N-t} [l_{(t)},b_{(m)}]z^{N-t-m-1} \\
    &= \sum_m \sum_{t=0}^N \binom{N}{t}(-1)^t [l_{(N-t)},b_{(m+t)}]z^{-m-1}
\end{split}
\]
The $m=-1$ coefficient is $-4a_{(N-2)} \neq 0$ for $N\ge 3$, forcing weak locality to fail. Since all the chosen generators are weakly local to one another, this contradictions Theorem 4.11(b) in \cite{Kac1999}.
\end{example}

\newpage

\section{Field algebras of non-symmetric lattices}  \label{sec: lattice}
We now build on the Heisenberg field algebra construction done in \Cref{sec: reconstruction} and develop a field algebra associated to any finite-rank lattice with a pairing. As in the above setup, let $L$ be a finite rank lattice with a (not necessarily symmetric) pairing $( \cdot, \cdot) : L \times L \rightarrow \mathbb Z$. Consider the complexification $\mathfrak h = L\otimes_{\mathbb Z} \mathbb C $.

We define for each $\lambda \in L$ a space $\pi_{\lambda}$ as in the symmetric case, defined by: \[
\begin{split}
h_n \cdot \ket \lambda = \delta_{n,0}\cdot (h,\lambda) \cdot \ket \lambda  \\
\Yop{h_n} \cdot \ket \lambda = \delta_{n,0}\cdot  (\lambda,h) \cdot \ket \lambda
\end{split}
\] for each $n \ge 0$. Each $\pi_\lambda$ is a module over the Fock space $\pi=\pi_0$, and its elements are generated by applying $h_{n}$ for $n < 0$, $h \in \frak{ h}$ to $\ket \lambda$. 

The translation operator acts by $T \cdot \ket \lambda = \lambda_{-1} \ket \lambda$. One can check that these extensions do not affect relations such as \Cref{lemma: heisenberg-op-commute}. Now let: $$V_L := \bigoplus_{\lambda \in L} \pi_\lambda $$

The reconstruction theorem (\Cref{thm : reconstruction}) will allow us to build a field algebra structure on $V_L$. Our candidate generators $S$ will be tuples of the form $(h^{(1)},\cdots, h^{(m)}, \lambda)$ where each $h^{(i)} \in \frak h$, $\lambda \in L$, and we have $a^{h^{(i)}}(z):= h^{(i)}(z)$ (resp $b^{h^{(i)}}(z)=(h^{(i)})^\text{op}(z)$), and $a^{\lambda}(z):= V_\lambda(z)$ (resp $b^{\lambda}(z):= \Yop{V}_\lambda(z)$) are to be defined.

So we only need to define $V_\lambda(z)$ and $\Yop{V}_\lambda(z)$ associated to $\ket \lambda $ for $\lambda \in L$. These fields cannot be random and should satisfy some formulas. Furthermore, since our vector space is spanned by elements of the form $\prod_i h^{\alpha_i}_{(j_{i})}\cdot \ket \lambda$ for $j_i < 0$, we only need that $(h(z),V_\lambda(z))$ is weakly local for all $\lambda \in L$ and $h \in \mathfrak h$. The reverse weak locality may not be true in general. 

In particular, one has that: \[
T \cdot \ket \lambda = \lambda_{-1}\ket \lambda \quad  \implies \quad \partial_zV_\lambda(z) = \normord{\lambda(z)V_\lambda(z)}
\]
OPE expansion (\Cref{cor: ope-on-source-fields}) should also give: \[
\begin{split}
    [\Yop{h_n},V_\lambda(z)] = (\lambda,h)z^nV_\lambda(z) \\
    [h_n,\Yop{V}_\lambda(z)] = (h,\lambda)z^n\Yop{V}_\lambda(z)
\end{split}
\]
We now give the candidates for our fields. We first define operators $S_\lambda,\Yop{S}_\lambda$ on $V_L$ as follows: For each $\nu,\lambda \in L$ we have that: 
\begin{equation}
    S_\lambda \cdot \ket \nu = c_{\lambda,\nu}\ket{\lambda+\nu}
\end{equation}
\begin{equation}
    \Yop{S}_\lambda \cdot \ket \nu = \Yop{c}_{\lambda,\nu} \ket{\lambda+\nu}
\end{equation}
where $c_{\lambda,\nu},\Yop{c}_{\lambda,\nu} \in \mathbb Z$ are coefficients that are to be determined. We then define for each $\lambda \in L$:

\begin{equation}
    \label{eq: exp-defn}
    V_\lambda(z) = S_\lambda \cdot z^{\lambda_0} \cdot \exp \bigg ( - \sum_{n < 0} \dfrac{\lambda_n}{n}z^{-n} \bigg ) \exp \bigg ( - \sum_{n > 0} \dfrac{\lambda_n}{n}z^{-n} \bigg )
\end{equation}
\begin{equation}
    \label{eq: exp-defn-op}
    \Yop{V}_\lambda(z) = \Yop{S}_\lambda \cdot z^{\Yop{\lambda}_0} \cdot \exp \bigg ( - \sum_{n < 0} \dfrac{\Yop{\lambda}_n}{n}z^{-n} \bigg ) \exp \bigg ( - \sum_{n > 0} \dfrac{\Yop{\lambda}_n}{n}z^{-n} \bigg ) 
\end{equation}
 where $z^{\lambda_0}$ acts on each $\pi_\nu$ by $z^{(\lambda,\nu)}$ and $z^{\Yop{\lambda}_0}$ acts on $\pi_\nu$ by $z^{(\nu,\lambda)}$. Both of these are well-defined operators since the second exponential eventually disappears when applied to a fixed vector in $V_L$. 
\theorem Let $L$ be a finite rank lattice with a bilinear form $(\cdot,\cdot) : L \times L \rightarrow \mathbb Z$. Suppose we pick $c_{\lambda,\nu}$ and $\Yop{c}_{\lambda,\nu}$ such that: \[
    \Yop{S}_\lambda S_\mu = (-1)^{(\mu,\lambda)}S_\mu\Yop{S}_\lambda 
    \] 
for any $\lambda,\mu \in L$. Then, the assignments \eqref{eq: exp-defn} and \eqref{eq: exp-defn-op} along with \eqref{eq: heisenberg-defn} extend to a field algebra structure on $V_L$ in the sense of \Cref{thm : reconstruction}. Furthermore, setting $c_{\lambda,\nu}=1$ and $c_{\lambda,\nu}^\text{op}=(-1)^{(\nu,\lambda)}$ satisfies the above condition, and we call the induced field algebra $(V_L,\ket 0, T,Y)$ the \textit{canonical field algebra structure associated to $L$}.
\begin{proof}
    One can see that $V_\lambda(z),\Yop{V}_\lambda(z)$ are fields for each $\lambda \in L$. We begin by showing some of the promised properties. Define: \[
    E_{1,\lambda}(z) = \exp \bigg ( - \sum_{n < 0} \dfrac{\lambda_n}{n}z^{-n} \bigg ) 
    \] and \[
    E_{2,\lambda}(z) = \exp \bigg ( - \sum_{n > 0} \dfrac{\lambda_n}{n}z^{-n} \bigg )
    \] and define $\Yop{E}_{1,\lambda}(z),\Yop{E}_{2,\lambda}(z)$ similarly. Then: \[
    \begin{split}
        \partial_z V_\lambda(z) &= S_\lambda \lambda_0 z^{\lambda_0-1}E_{1,\lambda}(z)E_{2,\lambda}(z) \\
        &+ S_\lambda z^{\lambda_0}\sum_{n < 0} \lambda_nz^{-n-1}E_{1,\lambda}(z)E_{2,\lambda}(z) \\
        &+ S_\lambda z^{\lambda_0}E_{1,\lambda}(z) \sum_{n > 0}\lambda_n z^{-n-1} E_{2,\lambda}(z) \\
        &= V_\lambda(z) \dfrac{\lambda_0}{z}+\lambda(z)_+V_\lambda(z)+V_\lambda(z)\lambda(z)_{-} - V_\lambda(z) \dfrac{\lambda_0}{z}\\
        &= \normord{\lambda(z)V_\lambda(z)}
    \end{split}
    \] Likewise $\partial_z\Yop{V}_\lambda(z) = \normord{\Yop{\lambda}(z)\Yop{V}_\lambda(z)}$. 

    Notice that a computation shows $[\Yop{h}_0,S_\lambda]=(\lambda,h)\cdot S_\lambda$ so that $[\Yop{h}_0,V_\lambda(z)] = (\lambda,h) V_\lambda(z)$ too. For $n > 0$, $\Yop{h}_n$ commutes with everything except with $E_{1,\lambda}(z)$, so if we write: \[
    E_{1,\lambda}(z) = \sum_{k \ge 0}\dfrac{N_\lambda(z)^k}{k!}
    \] for $N_\lambda(z) = -\sum_{m<0}\dfrac{\lambda_m}{m}z^{-m}$, so that $[\Yop{h}_n,N_\lambda(z)] = (\lambda,h)z^n$ and hence $[\Yop{h}_n,V_\lambda(z)] = (\lambda,h)z^nV_\lambda(z)$ for all $n\ge 0$. The other cases follow similarly, so that for all $n \in \mathbb{Z}$ we have: \[
    \begin{split}
        [\Yop{h}_n,V_\lambda(z)] &= (\lambda,h)z^nV_\lambda(z) \\
        [h_n,\Yop{V}_\lambda(z)]&=(h,\lambda)z^n\Yop{V}_\lambda(z)
    \end{split}
    \] In fact, by following the exact same format as above, we see ditto that: \[
    [h_n,V_\lambda(z)]=z^nV_\lambda(z) \cdot \big ( \textbf{1}_{n \ge 0}(h,\lambda)+\textbf{1}_{n<0}(\lambda,h) \big )
    \] and hence that $[h(z),V_\lambda(w)]_{z_-} = (h,\lambda)V_\lambda(w)\delta(z-w)_{z_-}$ so that $(h(z),V_\lambda(z))$ is weakly local for all $h \in \mathfrak h $ and $\lambda \in L$ by \Cref{lemma: weak-loc-formulation}. Likewise, $(\Yop{h}(z),\Yop{V}_\lambda(z))$ is weakly local, and $(h(z),\Yop{V}_\lambda(z))$ and $(\Yop{h}(z),V_\lambda(z))$ are local. 
    
    The last thing to show is that $(\Yop{V}_\lambda(z),V_\mu(z))$ is local. As in \cite{FrenkelBenZvi2004}, we see that: \[
    \bigg [  \sum_{n > 0} \dfrac{\Yop{\lambda}_n}{n}z^{-n},\sum_{m < 0} \dfrac{\mu_m}{m}w^{-m}  \bigg]
    = -\sum_{n > 0 } \dfrac{(\mu,\lambda)}{n}z^{-n}w^n = (\mu,\lambda)\log(1-\dfrac{w}{z})\] so by writing down $\Yop{V}_\lambda(z)V_\mu(w)$ and $V_\mu(w)\Yop{V}_\lambda(z)$ we see that for locality to hold it is enough that: \[
    \Yop{S}_\lambda S_\mu = (-1)^{(\mu,\lambda)}S_\mu\Yop{S}_\lambda 
    \] which holds either by assumption (in the first statement of the theorem) or by substitution for the canonical choice.
\end{proof} 

Notice that the field algebra we get may depend on our choice of $2$-cocycle $c_{\lambda,\mu},\Yop{c}_{\lambda,\mu}$. Another possible choice would have been to pick $c_{\lambda,\mu} = (-1)^{(\lambda,\mu)}$ and $\Yop{c}_{\lambda,\mu}=1$. For example, if $(a,b) = q(a,b)+q(b,a)$ is the symmetrization of some bilinear form $q : L \times L \rightarrow \mathbb Z$, the assignment $c_{a,b}=\Yop{c}_{a,b}=(-1)^{q(a,b)}$ gives a vertex algebra structure on $V_L$.  

\begin{example}
\label{ex: lattice}
Notice that it might not be the case that $(V_\lambda(z),h(z))$ is weakly local in a field algebra associated to a lattice. We continue with the same lattice as \Cref{ex: heisenberg}. We then have that: \[
[\ket{b}_{(m)},a_{(n)}] = (-1)^{\textbf{1}_{n < 0}} \ket{b}_{(m+n)}
\]
so that for all $N \ge 0$ we see that: \[
\begin{split}
    \res_z (z-w)^N [V_b(z),a(w)] &= \
\res_z\sum_{k=0}^N \sum_{m,n \in \mathbb{Z}} \binom{N}{k}(-1)^{N-k} [\ket b_{(n)}, a_{(m)}]  z^{k-n-1}w^{N-k-m-1} \\
    &= \sum_{m \in \mathbb{Z}} \ket b_{(N+m)} w^{-m-1} \cdot \sum_{k=0}^N \binom{N}{k} (-1)^{N-k} \cdot (-1)^{\textbf{1}_{m < k-N}} 
\end{split}
\] In particular, the coefficient at $m=-N$ is $2\cdot (-1)^N \cdot \ket b_{(0)} \neq 0$ since $\ket b_{(0)} a_{(-1)}\ket 0  =- \ket b_{(-1)} \ket 0 = - \ket b \neq 0$. 
\end{example}

\newpage 

\section{Commutant of field algebras}  \label{sec: centralizer}

Let us recall the definition of commutants of vertex algebras, see \cite[Theorem 5.1]{FrenkelZhu1992}. If $(V, \ket 0, T,Y)$ is a vertex algebra and $W$ is a subalgebra of $V$, the space: \[
W' := \{b \in V~|~ a_{(n)}b = 0 \text{ for all } a \in W\text{ and }n\ge0\} = \bigcap_{\substack{a \in W \\ n \ge 0}} \ker a_{(n)} 
\]
is also a subalgebra of $V$: If $b \in W'$, \Cref{prop : ope-goddard}(1) forces $[Y(a,z),Y(b,w)] = 0$, so that $b_{(m)}c\in W'$ for any $m \in \BZ$ and $c \in W'$. 

We now aim at defining an analog of the above picture for field algebras. Fix a vector space $V$ and consider a field algebra $(Y,Y^{\text{op}})$ on $V$ with vacuum vector $\ket{0}$ and translation operator $T$. Let $W$ be a subset of $V$. 

\begin{definition}[Strong commutant]
     Define the \textit{strong commutant} of $W$ in $V$ by: $$C^s_V(W) = C^s(W) := \big\{ b \in V ~|~ [a_{(n)},b_{(m)}] = 0 ~~ \forall a \in W ~;~\forall n \ge 0 ~; \forall m \in \mathbb{Z}\big\}$$ 
\end{definition}
We will show that $C^s(W)$ is a subfield algebra of $V$. The following proposition gives equivalent definitions for the strong commutant of a field algebra:
\begin{proposition}
    \label{prop: strong-commutant-equivalent-defs}
    Fix $a,b \in V$. The following are equivalent: \begin{enumerate}
        \item $[a_{(n)},b_{(m)}] = 0$ for all $n \ge 0$ and $m \in \mathbb{Z}$. 
        \item $Y(a,z)_{(n)}Y(b,z) = 0$ for all $n \ge 0$. 
        \item $(Y(a,z),Y(b,w))$ is weakly local and $a_{(n)}b = 0$ for all $n \ge 0$
    \end{enumerate}
\end{proposition}
\begin{proof} By \eqref{eq: nth-prod-coef} we have for all $n \ge 0$: \[
    Y(a,z)_{(n)}Y(b,z) = \sum_{m \in \mathbb{Z}} \sum_{k=0}^n \binom{n}{k}(-1)^{n-k}[a_{(k)},b_{(m-k)}] z^{n-m-1}
\]
This clearly gives $(1) ~\implies~ (2)$ and gives $(2) ~\implies~ (1)$ by an inductive argument on $n$. Both of them imply $(3)$ since $(2)$ implies the first part of $(3)$ by definition and $(1)$ implies the second part as $a_{(n)}b=[a_{(n)},b_{(-1)}]\ket 0 =0$ since $a_{(n)}\ket 0 =0$ for $n\ge 0$. Finally, $(3) \implies (2)$ since weak locality gives: \[Y(a,z)_{(n)}Y(b,z) = Y(a_{(n)}b,z) = 0\] for all $n \ge 0$ by \Cref{prop: topfields_nprod}.
\end{proof}
In the case of a vertex algebra $Y=Y^{\text{op}}$, the above conditions collapse to $a_{(n)}b = 0$ for all $n \ge 0$ as the $n$-th product axiom holds generally in that case \cite{BakalovKac2003}.
\begin{theorem}
    For any field algebra $(V,\ket{0},T,Y)$ and subset $W \subseteq V$, the strong commutant $C^s(W)$ forms a subfield algebra of $V$.
\end{theorem} 
\begin{proof} It is clear why $C^s(W)$ is a subvector space of $V$. It contains the vacuum $\ket{0}$ since $Y(\ket{0},z) = I_V$ commutes with every operator. It is closed under the translation operator $T$ as \[
[a_{(n)},(Tb)_{(m)}] = [a_{(n)},-mb_{(m-1)}] = -m[a_{(n)},b_{(m-1)}]=0
\] for all $a \in W, b \in C^s(W), n \ge 0$ and $m \in \mathbb{Z}$. We argue that it is closed under the action of $Y$. Fix $b,c \in C^ss(W)$ and $a \in W$. We will show from associativity that $[a_{(l)},(b_{(n)}c)_{(m)}] = 0$ for all $n,m \in \mathbb{Z}$ and $l \ge 0$. Fix a vector $v \in V$.  By \eqref{eq: associativity-coef}, we can find a large enough $N$ such that the following relation holds: \[
    \sum_{k=0}^N \binom{N}{k} (b_{(n+k)}c)_{(m-k)}x = \sum_{k=0}^N \sum_{j \ge 0} \binom{N}{k} \binom{-n-k+j-1}{j} b_{(n+k-j)} c_{(m+j-k)}x
\]
for $x \in \{v,a_{(l)}v\}$. Applying $a_{(l)}$ to the equation for $x=v$ and subtracting the equation for $x=a_{(l)}v$ gives: 
\[
\sum_{k=0}^N \binom{N}{k} [a_{(l)}, (b_{(n+k)}c)_{(m-k)}]v = 0
\] Here, $N$ depends on $b,c,v,a_{(l)}v$, so the relation holds for all $m,n$. But because $Y(b,z)$ is a field one may find $N_0$ for which $b_{(n)}c = 0$ for all $n > N_0$. The above equation at $n=N_0$ hence shows $[a_{(l)},(b_{(N_0)}c)_{(m)}]v=0$ for all $m \in \BZ$. At $n=N_0-1$, we get $[a_{(l)},(b_{(N_0-1)}c)_{(m)}]v+N \cdot [a_{(l)},(b_{(N_0)}c)_{(m-1)}]v = 0$, but the second term vanishes by the case $n=N_0$, giving us $[a_{(l)},(b_{(N_0-1)}c)_{(m)}]v=0$ for all $m \in \BZ$. In general, if we know that $[a_{(l)},(b_{(N')}c)_{(m)}]v=0$ for all $m \in \BZ$, taking $n=N'-1$ in the above equation gives vanishing of all higher order terms, forcing $[a_{(l)},(b_{(N'-1)}c)_{(m)}]v=0$ for all $m \in \BZ$. So $[a_{(l)},(b_{(n)}c)_{(m)}]=0$ for all $n,m \in \BZ$, as desired.
\end{proof} 

\begin{example}
    \label{ex: weak-com}
    If $(V,\ket 0, T,Y)$ is a vertex algebra and $W \subset V$, then $C^s(W) = C^s(\langle W\rangle )$, i.e we may replace $W$ by the sub-vertex algebra it generates in $V$: Indeed, it is clear that $C^s(W)\supset C^s(\langle W\rangle)$. On the other hand \[
    W \subset C^s(C^s(W)) \implies \langle W \rangle \subset C^s(C^s(W)) \implies C^s(\langle W\rangle) \supset C^s(C^s(C^s(W))) = C^s(W)
    \] where the last equality follows by \cite[Theorem 5.2]{FrenkelZhu1992}.
    
    This is not true for field algebras, i.e $C^s(W) \neq C^s(\langle W\rangle)$ in general. Consider for example the lattice $L = \mathbb Za \oplus \mathbb Zb \oplus \mathbb Z c$ with a form given by: \[
    Q = \begin{pmatrix}
        0 & 0 & 1\\
        1 & 0 & 0 \\
        0 & 0 & 0
    \end{pmatrix}
    \]
    Then, for the same reasons as in \Cref{ex: heisenberg}, one has that: \[
    Y(a_{(-1)}a_{(-1)}\ket 0 ,z) = \sum_{n \in \mathbb Z} \sum_{k+l=n-1} a_{(k)}a_{(l)}z^{-n-1}
    \] since $(a,a)=0$. Then, one has that $b \in C^s(a)$ but: \[
    \begin{split}
        [(a_{(-1)}a)_{(0)}, b_{(1)}] &= \sum_{k+l=-1} [a_{(k)}a_{(l)},b_{(1)}] \\
        &= -2 a_{(0)} \neq 0
    \end{split}
    \]
    since, for example, $a_{(0)} \ket c = \ket c \ne 0$. 
\end{example} 

In the above example, notice that $\ket a_{(n)} \ket c = 0$ for all $n \ge 0$ since $(a,c)=1\ge 0$. However: \[
[\ket a_{(0)}, \ket c_{(-1)}]\cdot c_{-1}\ket 0 = \ket{a+c} \neq 0
\] In particular, weak locality in condition (3) of \Cref{prop: strong-commutant-equivalent-defs} is not redundant. By dropping it, however, we get a weaker notion of commutants for exponentials $\ket \lambda$ in the context of lattice field algebras.

\subsection{Commutant for exponentials in lattice field algebras} Fix a lattice $L$ with a pairing $Q$ and let $V_L$ be the field algebra associated to this lattice as per \Cref{sec: lattice}. We will build a weaker notion of commutant for exponential terms $|\alpha \rangle$ in this field algebra. 

\begin{definition} (Commutant): 
    \label{def: weak-commutant}
    For each $\alpha \in L$, define the \textit{commutant} of $\alpha$ to be: \begin{align*} 
    C(\alpha) &:= \{v \in V_L~|~ |\alpha\rangle_{(n)}v = 0 \text{ for all } n \ge 0\}  \\
    &= \bigcap_{n \ge 0} \ker (|\alpha\rangle_{(n)})
    \end{align*}
\end{definition}

Naturally, $C^s(\alpha) \subseteq C(\alpha)$, and the inclusion can be strict as in \Cref{ex: weak-com}. Notice that $\ket \gamma \in C(\alpha)$ if and only if $(\alpha,\gamma) \ge 0$.

\begin{example}
    Fix a lattice $L$ with an arbitrary pairing and consider its associated field algebra. Consider $b = h_{-1}g_{-1}\cdot \ket \gamma$ for some $\gamma \in L$. Then for $n\ge 0$: \begin{align*}
        \ket \alpha_{(n)}b &= (\alpha,h)(\alpha,g)\cdot  \ket \alpha_{(n-2)} \ket \gamma \\
        &-(\alpha,h)g_{-1} \ket \alpha_{(n-1)}\ket \gamma - (\alpha,g) h_{-1} \ket \alpha_{(n-1)}\ket \gamma \\
        &+ h_{-1}g_{-1}\ket \alpha_{(n)}\ket \gamma
    \end{align*}
    For $b \in C(\alpha)$ to hold, one must necessarily have $\ket \gamma \in C(\alpha)$ since otherwise one can plug in $n = N_\gamma + 2,N_\gamma+1$ and $N_\gamma \ge 0$ successively to deduce a contradiction, where $N_\gamma$ is the largest integer for which $\ket \alpha_{N_\gamma} \ket \gamma \neq 0$. Notice furthermore that by our construction in \Cref{sec: lattice}, $\ket \alpha_{(r)}\ket \gamma = 0$ for any $r \ge -(\alpha,\gamma)$ and $\ket \alpha_{(-(\alpha,\gamma)-1)}\ket \gamma = \ket{\alpha+\gamma} \neq 0$. As a consequence of these facts: \begin{itemize}
        \item If $(\alpha,\gamma) \ge 2$, then $b \in C(\alpha)$ automatically.
        \item If $(\alpha,\gamma)=1$, then $b \in C(\alpha) \iff (\alpha,h)(\alpha,g)=0$.
        \item If $(\alpha,\gamma)=0$, then $b \in C(\alpha)\iff (\alpha,h)=(\alpha,g)=0$. 
    \end{itemize}
    
\end{example}

We now show that $C(\alpha)$ is a subfield algebra of $V_L$ for each $\alpha \in L$. In order to do so, we shall lift the condition on vectors to a condition on operators associated to said vectors. The classical approach in the context of vertex algebras is straightforward, since locality and OPE allow us to express $[Y(|\alpha \rangle,z),Y(v,w)]$ in terms of $Y(|\alpha\rangle_{(n)}v,w)$ as in \Cref{prop : ope-goddard}, and the latter terms vanish should $v$ be in the commutant of $\alpha$. We argue an analogous lifting in the nonsymmetric case. 

\begin{lemma}
    \label{lemma: wick}
    Let $\alpha \in L$ and $h \in \mathfrak h$, and consider $n,m \in \BZ$ such that either $n \ge 0$, or both $n,m < 0$. For any field $C(z)$, the following identity holds: \begin{align}
    Y(\ket \alpha , z)_{(n)}(\Yop{Y}(h_{-1} ,z)_{(m)}C(z)) &= -(\alpha,h)Y(\ket \alpha,z)_{(n+m)}C(z) \label{eq:wick-claim} \\
    &+ \Yop{Y}(h_{-1},z)_{(m)}(Y(\ket \alpha, z)_{(n)}C(z)) \notag
    \end{align}
\end{lemma}
\begin{proof}
    We separate cases based on the sign of $n$. If $n \ge 0$, this follows from Wick's theorem (Proposition 3.3(c) in \cite{Kac1998}) by noting that for $j \ge 0$: \[
    Y(\ket \alpha,z)_{(j)}\Yop{Y}(h_{-1}\ket 0,z) = -(\alpha,h)\cdot \delta_{j,0} \cdot Y(\ket \alpha,z)
    \]
    We now handle the case $n,m < 0$. Let $A(z) = Y(\ket \alpha,z)$, $B(z)=\Yop{Y}(h_{-1}\ket 0,z)$, and set $n'=-n-1 \ge 0$ and $m'=-m-1 \ge 0$. Then:
    \begin{align}
        (n'!)\cdot (m'!)\cdot A(z)_{(n)}(B(z)_{(m)}C(z)) &= \normord{\partial_z^{n'}A(z) \normord{\partial_z^{m'}B(z)C(z)}} \label{eq:wick-claim-proof}  \\
        &= \normord{\partial_z^{m'}B(z)\normord{\partial_z^{n'}A(z)C(z)}} \notag \\
        &+ C(z) [(\partial_z^{m'}B(z))_-, (\partial_z^{n'}A(z))_-] \notag \\
        &+[(\partial_z^{n'}A(z))_+, (\partial_z^{m'}B(z))_+]C(z) \notag
    \end{align}
    We now show the last two terms of \eqref{eq:wick-claim-proof} combine after rescaling to the first term in the right hand side of \eqref{eq:wick-claim}. Indeed: 
    \begin{align*}
        [(\partial_z^{n'}A(z))_+,(\partial_z^{m'}B(z))_+] &= \left [ \sum_{a < 0} \binom{-a-1}{n'}n'! \ket \alpha _{(a)}z^{-a-n'-1} ~,~ \sum_{b < 0} \binom{-b-1}{m'}m'! \Yop{h}_b z^{-b-m'-1} \right ] \\
        &= -(\alpha,h) \cdot m'!n'! \cdot \sum_{a,b<0} \binom{-a-1}{n'} \binom{-b-1}{m'} \ket \alpha_{(a+b)}z^{-(a+b)-m'-n'-2} \\
        &= -(\alpha,h)  \cdot m'! \cdot n'! \cdot \sum_{t < 0}\binom{-t-1}{m'+n'+1} \ket \alpha_{(t)}z^{-t-(m'+n'+1)-1}
    \end{align*}
    where the last summation is precisely the holomorphic part of $\partial_z^{-m-n-1}A(z) \cdot \frac{1}{(-m-n-1)!}$ thanks to the Chu-Vandermonde identity. The other term is computed analogously and yields: \[
    [(\partial_z^{m'}B(z))_-,(\partial_z^{n'}A(z))_-] = -(\alpha,h) \cdot m'! \cdot n'! \cdot \sum_{t \ge 0} \binom{-t-1}{m'+n'+1} \ket \alpha_{(t)}z^{-t-(m'+n'+1)-1}
    \] which is precisely the singular part of $\partial_z^{-n-m-1}A(z)\cdot \frac{1}{(-m-n-1)!}$
\end{proof}

\begin{remark}
    \label{remark: wick-heisenberg}
    A similar argument yields that for any $h,g \in \mathfrak h$ and integers $m,n<0$, one has that: \[
    \Yop{Y}(h_{-1},z)_{(n)}(\Yop{Y}(g_{-1},z)_{(m)}C(z)) = \Yop{Y}(g_{-1},z)_{(m)}(\Yop{Y}(h_{-1},z)_{(n)}C(z))
    \] for any field $C(z)$.
\end{remark}

The significance of \Cref{lemma: wick} is that it allows us to mimic the computation of $n$-th products in terms of the computations in $V_L$ of their associated vectors. 

\begin{proposition}
    \label{prop: finishing-prop}
    For each $\alpha \in L$ and $v \in C(\alpha)$, one has for all $n \ge 0$: \[
    Y(|\alpha \rangle, w)_{(n)}\Yop{Y}(v,w)= 0
    \] and hence by locality: \[
    [Y(|\alpha \rangle, z), \Yop{Y}(v,w)] = 0
    \]
\end{proposition}
\begin{proof}
    Consider the map $\sigma : \pi_0 \rightarrow \pi_0[x]$ given on generators $h_{m}~,~ m < 0$ for $h \in \frak h$ by: \[
    \sigma(\ket 0) = \ket 0 \quad ; \quad \sigma(h_{m}) = h_{m}-(\alpha,h)x^{-m}
    \] and extended to $\pi_0$ multiplicatively. For any $w \in \pi_0$ and $\ell \in \BZ$, let $m_\ell(w) \in \pi_0$ denote the $x^\ell$'th coefficient of $\sigma(w)$. First, we claim that for any $\gamma \in L,w\in \pi_0$: \begin{equation}
    \label{eq: commtrick}
    \ket \alpha_{(n)} (w\cdot \ket \gamma) = \sum_{\ell \in \BZ} m_\ell(w) \cdot \ket \alpha_{(n-\ell)} \ket \gamma 
    \end{equation} The right hand side is always well-defined since $m_\ell(w)$ vanishes for $|\ell|$ large enough. Since \eqref{eq: commtrick} is additive in $w$, it's enough to show it for $w = h_{i_1}^{(j_1)}\cdots h_{i_k}^{(j_k)}$ with each $h^{(j_q)} \in \frak h$, which we do by induction on $k$, with $k=0$ following by definition of $\sigma$. For $w = h_{i}w'$ with $i < 0$ and $h \in \frak h$, we have: \begin{align*}
    \ket \alpha_{(n)}(w\cdot \ket \gamma) &= -(\alpha,h) \ket \alpha_{(n+i)}w'\cdot \ket \gamma + h_i\cdot (\ket \alpha_{(n)}w'\cdot \ket \gamma) \\
    &=-(\alpha,h) \sum_{\ell \in \BZ} m_\ell(w') \ket \alpha_{(n+i-\ell)} \ket \gamma + h_i \sum_{\ell \in \BZ} m_\ell(w')\ket \alpha_{(n-\ell)} \ket \gamma \\
    &= \sum_{\ell \in \BZ}\big ( h_im_\ell(w')-(\alpha,h)m_{\ell +i}(w')\big ) \ket \alpha_{(n-\ell)}\ket \gamma\\
    &= \sum_{\ell \in \BZ} m_{\ell}(w) \ket \alpha_{(n-\ell)} \ket \gamma 
    \end{align*} as desired. 

    Let $glf(V_L)$ denote the space of fields on $V_L$ (see \cite{BakalovKac2003}). By \Cref{remark: wick-heisenberg}, we have a map: \[
    \varphi : \pi_0 \otimes glf(V_L) \rightarrow glf(V_L)
    \] given by: \[h^{(i_1)}_{-j_1}h^{(i_2)}_{-j_2}\cdots h^{(i_k)}_{-j_k}\otimes C(z) \mapsto \Yop{Y}(h^{(i_1)}_{-1}\ket 0,z)_{(-j_1)}(\Yop{Y}(h_{-1}^{(i_2)}\ket 0,z)_{(-j_2)}(\cdots C(z))\cdots )\]
    and extended linearly. We inductively claim that for any $w \in \pi_0$ and field $C(z)$ we have: \[
    Y(\ket \alpha,z)_{(n)} \varphi(w \otimes C(z)) = \sum_{l \in \BZ} \varphi\big (m_\ell(w) \otimes Y(\ket \alpha,z)_{(n-\ell)}C(z) \big )
    \] for similar reasons as above: It clearly holds for $w = \ket 0$, and assuming $w = h_i w'$, then by \Cref{lemma: wick}: \begin{align*}
        Y(\ket \alpha,z)_{(n)}\varphi(w \otimes C(z)) &= Y(\ket \alpha,z)_{(n)}(\Yop{Y}(h_{-1},z)_{(i)}\varphi(w'\otimes C(z))) \\
        &=-(\alpha,h)Y(\ket \alpha,z)_{(n+i)}\varphi(w'\otimes C(z)) \\
        &+\Yop{Y}({h_{-1}},z)_{(i)}(Y(\ket \alpha,z)_{(n)}\varphi(w'\otimes C(z))) \\
        &= -(\alpha,h)\sum_{\ell \in \BZ}\varphi(m_\ell(w')\otimes Y(\ket \alpha,z)_{(n-\ell+i)}C(z)) \\
        &+ \sum_{\ell \in \BZ}\Yop{Y}(h_{-1},z)_{(i)}\big (\varphi(m_\ell(w') \otimes Y(\ket \alpha,z)_{(n-\ell)}C(z)) \big ) \\
        &= \sum_{\ell \in \BZ} \varphi\big(m_\ell(w) \otimes Y(\ket \alpha,z)_{(n-\ell)}C(z)\big )
    \end{align*} as desired. \\

    In particular, if $v \in C(\alpha)$ is chosen arbitrarily, then since $\ket \alpha_{(n)}$ sends $\pi_\gamma$ to $\pi_{\alpha+\gamma}$ for each $\gamma \in L$, we may assume that $v=w \ket \gamma$ for some Heisenberg $w$. Then one can expand:
    \begin{equation}
        \label{eq: expansion}
        \ket \alpha_{(n)}v = \sum_{\ell \ge 0}^N m_\ell \cdot \ket \alpha_{(n-\ell)} \ket \gamma = 0 \quad ; \quad m_\ell =  m_\ell(w) \quad \text{ for some } N >  0
    \end{equation}
    By definition (the reconstruction theorem), we have that $\Yop{Y}(v,z)=\varphi(w \otimes \Yop{Y}(\ket{\gamma}, z))$, hence by the above we have: \[
    Y(\ket \alpha,z)_{(n)}\Yop{Y}(v,z) = \sum_{\ell \ge 0}^N \varphi \left (m_\ell \otimes (Y(\ket \alpha,z)_{(n-\ell)}\Yop{Y}(\ket \gamma,z)) \right )
    \]
    Notice also that $Y(\ket \alpha,z)_{(n)}\Yop{Y}(\ket \gamma,z) = 0$ for $n \ge -(\alpha,\gamma)$ since: \[
    Y(\ket \alpha,z)\Yop{Y}(\ket \gamma,w) = (z-w)^{(\alpha,\gamma)}f(z,w)
    \] where: \[
    f(z,w) = S_\alpha \Yop{S}_\gamma  z^{\alpha_0}w^{\Yop{\gamma}_0} \exp\left (-\sum_{n < 0}  \dfrac{\alpha_nz^{-n}+\Yop{\gamma}_nw^{-n}}{n} \right )\exp\left (-\sum_{n > 0} \dfrac{\alpha_nz^{-n}+\Yop{\gamma}_nw^{-n}}{n} \right )
    \] is regular in $z-w$. In particular, if $\ell \le (\alpha,\gamma)$ we have: \[
    Y(\ket \alpha,z)_{(n-\ell)}\Yop{Y}(\ket \gamma,z) =0 ~\implies ~ \varphi(m_\ell \otimes Y(\ket \alpha,z)_{(n-\ell)}\Yop{Y}(\ket \gamma,z)) = 0
    \] where the first equality follows from $n-\ell\ge -(\alpha,\gamma)$. \\

    Now suppose that $\ell > (\alpha,\gamma)$. We claim that $m_\ell=0$ since $v \in C(\alpha)$, which would imply the result. Indeed, if not, truncate \eqref{eq: expansion} in such a way that $m_N \neq 0$ and $N > (\alpha,\gamma)$ is maximal with $m_N \neq 0$. Notice that: \[
    \ket \alpha_{(N-(\alpha,\gamma)-1)} v = \sum_{\ell =0}^N m_\ell \cdot \ket \alpha_{(N-(\alpha,\gamma)-\ell-1)}\ket \gamma = m_N \ket \alpha_{(-(\alpha,\gamma)-1)}\ket \gamma = m_N\cdot \ket{\alpha+\gamma} = 0
    \] where the first $N$ terms vanish since $\ket \alpha_{(n)}\ket \gamma = 0$ for $n \ge -(\alpha,\gamma)$. Thus $m_N=0$, contradiction.
\end{proof}

\begin{theorem}
    \label{thm: weak-commutant}
    Let $L$ be a finite rank lattice with bilinear form $(\cdot,\cdot) : L \times L \rightarrow \BZ$, and let $(V_L,\ket 0, T,Y)$ be the (canonical) field algebra structure associated to $L$. Then for each $\alpha \in L$, the subspace $C(\alpha)$ is a subfield algebra of $V_L$. 
\end{theorem}
\begin{proof}
    $C(\alpha)$ is clearly a linear subspace of $V_L$. State-field correspondence axioms imply that $C(\alpha)$ contains $\ket 0$ and is closed under $T$. Finally, \Cref{prop: finishing-prop} implies that $C(\alpha)$ is closed under $\Yop{Y}$: For $v \in C(\alpha)$, we have: \[
    [Y(\ket \alpha,z), \Yop{Y}(v,w)] = 0
    \] hence for any $v' \in C(\alpha), n \ge 0, m\in \BZ$ we have: \[
    \ket \alpha_{(n)}(\Yop{v}_{(m)}v') = \Yop{v}_{(m)}(\ket \alpha_{(n)}v') = 0
    \] so $\Yop{v}_{(m)}v' \in C(\alpha)$. By \Cref{prop: s2fprop} (4), we see that $C(\alpha)$ is also closed under $Y$.
\end{proof}

It is not true in general that $\bigcap_{ n \ge 0} \ker(v_{(n)})$ is a subfield algebra of a field algebra $V$ for any $v \in V$, as the next example showcases. 

\begin{example}
    \label{ex: comm-not-subfield}
    Let $U$ be the standard Heisenberg vertex algebra structure associated to $L=\BZ$ with pairing given by integer products with generator $\alpha$. Let $\sigma$ be the order $2$ field algebra automorphism of $U \otimes_\BC U$ swapping the two coordinates, and consider: \[
    V = (U \otimes U)  \sharp \langle 1,\sigma\rangle 
    \] (See Example 3.6 in \cite{BakalovKac2003}). Let $a = (\alpha_{-1} \otimes \ket 0) \otimes 1$, $b=(\ket 0 \otimes \ket 0) \otimes \sigma$ and $c = (\ket 0 \otimes \alpha_{-1}) \otimes 1$. Then: \[
    Y(a,z)b = Y(\alpha_{-1} \otimes \ket 0, z)(\ket 0 \otimes \ket 0) \otimes \sigma = (Y(\alpha_{-1},z) \ket 0 \otimes Y(\ket 0,z) \ket 0) \otimes \sigma
    \] Likewise $Y(a,z)c$ is also holomorphic. However: \[
    Y(b,z)c = Y(\ket 0 \otimes \ket 0 , z)(\alpha_{-1} \otimes \ket 0) \otimes \sigma
    \] and hence $b_{(-1)}c=( \alpha_{-1} \otimes \ket 0) \otimes \sigma$. In particular, notice that: \[
    Y(a,z)b_{(-1)}c = Y(\alpha_{-1}\otimes \ket 0,z)(\alpha_{-1} \otimes \ket 0) \otimes \sigma = (Y(\alpha_{-1},z)\alpha_{-1} \otimes \ket 0) \otimes \sigma  
    \] is not holomorphic. 
\end{example}

\newpage

\section{Lattice field algebra associated to quiver representations}

We now discuss the motivating example of this paper. Theorem 2.6 of \cite{GrossJoyceTanaka2020} details a general approach to associate a graded vertex algebra to the homology of moduli stacks $\CM$ associated to certain $\BC$-linear abelian or triangulated categories, so long as one provides the extra data of a monoidal action $\Phi : \CM \times \CM \rightarrow \CM$, a $B\BG_m$-action and a perfect complex $\Theta \in \Perf(\CM\times \CM)$ satisfying some compatibility conditions. 

For example, if $\overline{\CM}$ is the moduli stack of complexes of representations of a quiver, Theorem 5.19 in \cite{Joyce2019} shows one can endow $H_*(\overline{\CM})$ with the structure of a lattice vertex algebra. The choice of perfect complex here is given by: \[
\Theta = (\overline{\CE} xt)^\vee \oplus \sigma^*(\overline{\CE} xt)
\]  where $\sigma : \overline{\CM} \times \overline{\CM } \rightarrow \overline{\CM} \times \overline{\CM}$ is the automorphism swapping coordinates and $\overline{\CE} xt$ is a natural choice of perfect complex described below. The symmetrization is required since a (super)vertex algebra structure can only be associated to lattices with symmetric pairings (see Section 5 in \cite{FrenkelBenZvi2004}). 

One can drop the symmetrization above to get a braided vertex algebra, see \cite{Latyntsev2021}. In this section, we view this as a lattice field algebra structure on the homology of the moduli stacks of complexes of quiver representations. We begin by recalling some basic definitions. 

\begin{definition}
    \label{def: quiver} (Quiver): 
    A \textit{quiver} is a quadruple $Q=(Q_0,Q_1,s,t)$ where $Q_0$ is a finite set of vertices, $Q_1$ is a finite set of directed edges between vertices in $Q_0$ and $s,t : Q_1 \rightarrow Q_0$ specify the source and target of each edge. 
    
    A \textit{$\BC$-representation of} $Q$ is a collection $\{V_v, \varphi_e\}_{v\in Q_0,e\in Q_1}$ of finite-dimensional $\BC$-vector spaces $V_v$ associated to each vertex $v \in Q_0$ along with linear maps $\varphi_e : V_{s(e)}\rightarrow V_{t(e)}$ for each $e \in Q_1$.  
\end{definition} 

Quiver representations naturally form an abelian category $\modCQ{Q}$ where a morphism $\rho$ from $\{V_v,\varphi_e\}_{v\in Q_0,e\in Q_1}$ to $\{W_v, \psi_e\}_{v\in Q_0, e\in Q_1}$ is defined as a collection of linear maps $\rho_v : V_v \rightarrow W_v$ such that for each $e \in Q_1$ the following diagram commutes: \[
\begin{tikzcd}
    V_{s(e)} \arrow[d, "\rho_{s(e)}"'] \arrow[r, "\varphi_e"] & V_{t(e)} \arrow[d, "\rho_{t(e)}"] \\
    W_{s(e)} \arrow[r, "\psi_e"'] & W_{t(e)}
\end{tikzcd}
\]

We can then form the (triangulated) bounded derived category $D^b \modCQ{Q}$ corresponding to complexes of quiver representations. Let $\CM$ be the Artin $\BC$-stack with $\BC$-points corresponding to isomorphism classes of quiver representations. Likewise, let $\overline{\CM}$ denote the higher $\BC$-stack with $\BC$-points corresponding to objects in the triangulated category $D^b\modCQ{Q}$. We will be interested in the latter from now on.

We begin by recalling the homology of $\overline{\CM}$ which we take as the homology of the underlying topological realization $\overline{\CM}^{\text{top}}$ with $\BC$-coefficients. If $K(Q)$ denotes the Grothendieck group of the bounded derived category $D^b\modCQ{Q}$, then there is an isomorphism: \[
\textbf{dim} : K(Q) \xrightarrow{\cong} \BZ^{|Q_0|}
\] defined on a complex by sending at each vertex its Euler characteristic. Notice furthermore that the induced map: \[
\overline{\CM}(\BC) \twoheadrightarrow \BZ^{|Q_0|}
\] is also surjective and locally constant, which hence gives a decomposition of $\overline{\CM}$ into components: \[
\overline{\CM} = \coprod_{d \in \BZ^{|Q_0|}} \overline{\CM}_d
\]
Then, \cite{Joyce2019} show in Proposition 5.13 that there is a homotopy equivalence: \[
\overline{\CM}_d \simeq \prod_{v \in Q_0} \Perf_\BC^{d(v)}
\] (where $\Perf_\BC^r$ denotes perfect complexes of rank $r$) as well as an identification: \[
H_*(\overline{\CM}_d;\BC) \cong \BC[b_{d,v,i}~,~ i=1,2,\ldots ~,~ v \in Q_0]
\] with each $b_{d,v,i}$ of degree $2i$, and hence identifies $H_*(\overline{\CM};\BC)$ with the underlying vector space of a lattice field algebra associated to $\BZ^{|Q_0|}$. 

Furthermore, the quotient $\BZ^{|Q_0|}$ has a natural pairing known as the \textit{Euler pairing} which stems from the underlying quiver and is given by: 
\[
\chi(\alpha,\beta) = \sum_{v \in Q_0}  \alpha(v) \beta(v)  - \sum_{e \in Q_1} \alpha(s(e))\cdot \beta(t(e))
\] where $\alpha(v)$ denotes the $v$'th coordinate of $\alpha \in \BZ^{|Q_0|}$. We denote by $L_Q$ the lattice $\BZ^{|Q_0|}$ with the Euler pairing. 

We now proceed as follows. First, we explain that $H_*(\overline{\CM};\BC)$ has a structure of a field algebra by specifying a monoidal structure $\overline{\Phi}$, a $B\BG_m$ action and a perfect complex $\Theta$. This data will be enough to endow $H_*(\overline{\CM};\BC)$ with the structure of a braided vertex algebra, and hence a field algebra, as per \cite{Latyntsev2021} Section 3.7. 

Second, we will argue that the above identification is an isomorphism of field algebras between the field algebra given on $H_*(\overline{\CM};\BC)$ and the lattice field algebra associated with $L_Q$. 

We follow Section 5.3 from \cite{Joyce2019} for the first part. 

\subsection{Field algebra structure on $H_*(\overline{\CM};\BC)$} We describe the field algebra structure on $H_*(\overline{\CM};\BC)$ by describing the necessary data for said structure. We remark that the quotient of the Grothendieck group was described above using $\textbf{dim}$, and that the choice of signs is irrelevant in our case since we are aiming for a field algebra. This justifies why signs are not present on $Y$ for field algebras associated to lattices in \Cref{sec: lattice} and were instead pushed to $\Yop{Y}$.

We begin by describing the perfect complex $\Theta \in \Perf(\overline{\CM} \times \overline{\CM})$. In what follows, we pick $\Theta = \overline{\CE}xt^\vee$ where $\overline{\CE}xt$ is the natural complex whose cohomology at $\BC$-points of $\overline{\CM} \times \overline{\CM}$ recovers Ext groups in the triangulated category $D^b(\modCQ{Q})$: \[
H^i(\overline{\CE}xt|_{E,F}) \cong \Ext^i(E,F) := \Hom_{D^b(\modCQ{Q})}(E,F[i]) 
\] for $E,F \in D^b(\modCQ{Q})$. Furthermore, we see by equation (5.30) in \cite{Joyce2019} that: \[
\rank \overline{\CE}xt_{d,e} = \chi(d,e)
\] for $d,e \in \BZ^{|Q_0|}$, where $\overline{\CE}xt_{d,e} = \overline{\CE}xt|_{\overline{\CM}_{d} \times \overline{\CM}_e}$. 

The monoidal structure is given by a map $\overline{\Phi} : \overline{\CM} \times \overline{\CM} \rightarrow \overline{\CM}$ which restricts on each $\overline{\CM_d} \times \overline{\CM_e}$ component to maps $\overline{\Phi}_{d,e} : \overline{\CM}_d \times \overline{\CM}_e \rightarrow \overline{\CM}_{d+e}$, and the $B\BG_m$-action is given by a map $\overline{\Psi} : [\star / \BG_m] \times \overline{\CM} \rightarrow \overline{\CM} $ restricting to $\overline{\Psi}_d : [\star/\BG_m] \times \overline{\CM}_d \rightarrow \overline{\CM}_d$ on each component. 

We now explain what these maps are. On each factor, we have maps: \[
 \Perf_\BC^r \times \Perf_\BC^s \rightarrow \Perf_{\BC}^{r+s} \quad \quad ; \quad \quad (\mathscr
F^\bullet,\mathscr G^\bullet) \mapsto \mathscr F^\bullet \oplus \mathscr G^\bullet
\] which hence extend to a map: \[
\prod_{v \in Q_0} \Perf_\BC^{d(v)} \times \Perf_{\BC}^{e(v)} \rightarrow \prod_{v \in Q_0} \Perf_\BC^{(d+e)(v)}
\] which defines $\overline{\Phi}_{d,e}$. 

Likewise, the map $\overline{\Psi_d}$ is induced on each factor by: \[
[\star/\BG_m] \times \Perf_\BC^{d(v)} \rightarrow \Perf_\BC^{d(v)} \quad \quad ; \quad \quad (L,\mathscr F^\bullet) \mapsto L \otimes \mathscr F^\bullet 
\]

As per \cite{Latyntsev2021}, this data determines a field algebra structure on $H_*(\overline{\CM};\BC)$ as follows: \begin{itemize}
    \item The vacuum is given by the zero object $\ket 0 := [0] \in H_0(\overline{\CM}_0;\BC)$ and in particular lives in the homology of $\overline{\CM}_0$. 
    \item The translation operator $T$ acts on homology by: \[
    T(u) := \overline{\Psi}_*(t \times u) \in H_*(\overline{\CM};\BC)
    \] where $t\in H_2([\star/\BG_m])$ is the dual of the first Chern class.
    \item For homology classes $u \in H_*(\overline{\CM}_d;\BC)$ and $v \in H_*(\overline{\CM}_e;\BC)$ respectively, the state-field correspondence is given as follows (compare with equation (3.33) in \cite{Joyce2019}): \[
    Y(u,z)v := z^{\text{rank } \Theta_{d,e}} \cdot \overline{\Phi}_* \left ( (e^{zT} \otimes \id ) \cdot c_{z^{-1}}(\Theta) \cap ( u \times v)\right )
    \] where: \[
    c_{z^{-1}} := \sum_{i \ge 0} z^{-i}c_i
    \] is the total Chern class operator. Notice in particular $\text{rank } \overline{\CE}xt_{d,e} = \chi(d,e)$ due to our choice of perfect complex.
\end{itemize}

\begin{remark}
    The convention in \cite{Latyntsev2021} is for the complex $\Theta$ to have weights $(-1,1)$ with respect to the $B\BG_m\times B\BG_m$-action on $\overline{\CM}\times \overline{\CM}$, opposite to $(1,-1)$ which both \cite{Joyce2019} and we use. 
\end{remark}

We now clarify this structure. We begin by making our identification more specific. Since for each $d \in L_Q$ we have that $\overline{\CM}_d \simeq \prod_{v \in Q_0}\Perf_\BC^{d(v)}$, then if $\mathscr V_v$ denotes the universal perfect complex on each term in the product, we have an isomorphism of $\BC$-algebras: \[
H^*(\overline{\CM}_d;\BC) \cong \BC[\beta_{d,v,i} ~,~ v \in Q_0~,~ i \ge 1]
\] where each $\beta_{d,v,i} = \text{ch}_i(\mathscr V_v)$ are the Chern characters, see equation (5.33) in \cite{Joyce2019}. We can also write down an isomorphism of vector spaces: \[
H_*(\overline{\CM}_d;\BC) \cong \BC[b_{d,v,i}~,~ v\in Q_0~,~ i\ge 1]
\] where the monomials $b_{d,v,i}$ are duals to the Chern characters $\beta_{d,v,i}$, see Proposition 5.13 and equation (5.34) in \cite{Joyce2019}.
\begin{lemma}
    \label{lemma: chern-characters}
    The Chern characters of $\Theta=\overline{\CE}xt^\vee$ are given on $\overline{\CM_d} \times \overline{\CM_e}$ by: \[
    \text{ch}_i(\Theta_{d,e}) = \sum_{v,w \in Q_0} \sum_{\substack{j,k \ge 0 \\ j+k = i}} (-1)^k \chi(\textbf{1}_v,\textbf{1}_w) \cdot \beta_{d,v,j} \boxtimes \beta_{e,w,k}
    \] where $\textbf{1}_v \in \BZ^{Q_0}$ denotes the vector with entry $1$ at $v$ and $0$ everywhere else, and hwere $\beta_{d,v,0}$ acts by scaling by $d(v)$. 
\end{lemma}
\begin{proof}
    Follows by dualizing then taking Chern characters of equation (5.54) as in Proposition 5.15 of \cite{Joyce2019}. 
\end{proof}

\subsection{Lattice field algebra isomorphism} We are now ready to show that the field algebra structure on $H_*(\overline{\CM};\BC)$ is isomorphic to $V_{L_Q}$ for each quiver $Q$. 

Recall that each summand $\pi_d$ of $V_{L_Q}$ is generated as a $\BC$-vector space by elements of the form \[
v^{(j_1)}_{i_1} \cdot v_{i_2}^{(j_2)} \cdot \cdots \cdot v_{i_n}^{(j_n)} \cdot \ket d
\] where $v^{(j_\ell)} \in L_Q$ with $i_\ell < 0$. We denote for each $v \in Q_0$ and $i \ge 1$ the Heisenberg operator ${(\textbf{1}_v)}_{-i}$ on $\pi_d$ by $v_{i,d}$  
\begin{theorem}
    \label{thm: quiv-hom-is-lattice-field-algebra}
    For any finite quiver $Q$, there is an isomorphism of field algebras: \[
    \xi : H_*(\overline{\CM};\BC) \xrightarrow{\cong } V_{L_Q}
    \] given on each $d \in L_Q$ by linear maps: \[
    \xi_d : H_*(\overline{\CM}_d;\BC) \rightarrow \pi_d
    \] determined by: \[
    b_{d,v_1,i_1}\cdot b_{d,v_2,i_2}\cdot \cdots \cdot b_{d,v_n,i_n} \mapsto  (v_1)_{i_1,d}\cdot (v_2)_{i_2,d}\cdot \cdots \cdot (v_n)_{i_n,d} \cdot \ket d
    \] and extended linearly.
\end{theorem}
\begin{proof}
    The map is clearly bijective, so we only have to show that $\xi$ respects the field algebra structure. The vacuum $[0] \in H_0(\overline{\CM_0};\BC)$ is sent to $\ket 0$. So we only have to check compatibility, which boils down to computing the state-field correspondence $Y$ for $H_*(\overline{\CM};\BC)$. \\

    Since \Cref{lemma: chern-characters} holds, we can repeat the argument of Theorem 5.19 in \cite{Joyce2019} since the calculations of section 5.3.5 hold (with $a_{vw}$ replaced by $\chi(\textbf{1}_v,\textbf{1}_w)$). Thus we can identify the generators of $H_*(\overline{\CM};\BC)$ with those of $V_{L_Q}$, at which point we are done by \Cref{thm : reconstruction}.
\end{proof}

\begin{remark}
    One can make $\xi$ into a graded isomorphism by changing the grading on the homology $H_*(\overline{\CM};\BC)$ as done in \cite{Joyce2019}. 
\end{remark}

If we fix a dimension vector $d \in \BZ^{|Q_0|}$ (really an ``Euler characteristic vector''), we can interpret the commutant $C(d)$ (as defined in \Cref{def: weak-commutant}) in homological language by using Chern classes of $\Theta$. 

\begin{proposition}
    \label{prop: commutant=chern-class-capping}
    If $Q$ is a finite quiver and $d \in \BZ^{|Q_0|}$, then the bijection $\xi_e^{-1}$ in \Cref{thm: quiv-hom-is-lattice-field-algebra} sends $C(d) \cap \pi_e$ to: \[
    \{u \in H_*(\overline{\CM}_e;\BC) ~;~  c_i(\Theta_{d,e}) \cap  (\pt_d \times u) = 0 ~ \text{ for all } i > \chi(d,e) \}
    \] where $\pt_d \in H_0(\overline{\CM}_d;\BC)$ is a generator.
\end{proposition}
\begin{proof}
    An element $\psi \in \pi_e$ is in $C(d)$ if and only if: \[
    z^{\chi(d,e)}\exp(-\sum_{n > 0}\dfrac{d_n}{n}z^{-n})\psi
    \] is holomorphic in $z$. By the identification $\xi_{e}^{-1}$, we see that $d_n$ acts on $H_*(\overline{\CM}_e;\BC)$ by capping with $\sum_{v \in Q_0} \chi(d,\textbf{1}_v) n!\beta_{e,v,n}$. Thus $\xi_e^{-1}(\psi)=\eta \in \xi_e^{-1}(\pi_e \cap C(d))$ if and only if: \begin{equation} \label{eq:holo}
    z^{\chi(d,e)} \exp \left (-\sum_{n>0}\sum_{v \in Q_0} (n-1)!\cdot \chi(d,\textbf{1}_v)\cdot  \beta_{e,v,n}   \cdot z^{-n} \right )  \cap \eta 
    \end{equation} is holomorphic. Pick a generator $\pt_d \in H_0(\overline{\CM}_d;\BC)$. Then for each $j \ge 1$, $\beta_{d,v,j} \cap \pt_d=0$ and $\beta_{d,v,0} \cap (-)$ acts by scaling by $d(v)$. By identifying $\pt_d \times \eta$ with $\eta$, we can then invoke \Cref{lemma: chern-characters} to see that \eqref{eq:holo} equals:\[
    z^{\chi(d,e)} \exp \left ( \sum_{j \ge 1} (-1)^{j-1} \cdot(j-1)!\cdot  z^{-j}\cdot \ch_j(\Theta_{d,e})\right ) \cap (\pt_d \times \eta)
    \]  But one has: \[
    \sum_{j \ge 1} (-1)^{j-1} \cdot(j-1)!\cdot  z^{-j}\cdot \ch_j(\Theta_{d,e})  = \log c_{z^{-1}}(\Theta_{d,e}) = \log \left (  \sum_{i \ge 0} z^{-i} c_i(\Theta_{d,e}) \right ) 
    \] from which the result follows. 
\end{proof}

Let $v \in Q_0$ and suppose one can find $d_v \in \BZ^{|Q_0|}$ for which $\chi(d_v,\textbf{1}_w)=\delta_{vw}$, so that $\chi(d_v,e)=e(v)$ is the $v$-th coefficient of $e \in \BZ^{|Q_0|}$. Then \eqref{eq:holo} becomes: \begin{align*}
    z^{e(v)}\exp \bigg ( -\sum_{n > 0} (n-1)! \cdot \beta_{e,v,n}\cdot z^{-n}\bigg )  &= z^{e(v)}\exp(\sum_{n > 0}(-1)^{n-1} \cdot (n-1)! \cdot z^{-n}\cdot \ch_n(\mathscr V_v^\vee) )  \\
    &= \sum_{n \ge 0} z^{e(v)-n} c_n(\mathscr V_v^\vee) \\
    &= \sum_{n \ge 0}z^{e(v)-n}(-1)^{n}c_n(\mathscr V_v)
\end{align*}

The condition $\eta \in \xi_{e}^{-1}(C(d_{v}) \cap \pi_{e})$ then collapses to $c_i(\mathscr V_{v})\cap \eta = 0$ for $i> \chi(d_v,e)= e(v)$, where $\mathscr V_{v}$ is the universal perfect complex of $\Perf_{\BC}^{e(v)}$. In particular, we deduce the following:

\begin{corollary}
    \label{cor: moduli-stack-of-reps-is-commutant}
    Let $Q$ be a finite quiver with no oriented cycles. Let $\CM$ be the Artin $\BC$-stack of quiver representations associated to $Q$. Then, for each $v \in Q_0$, one can find $d_v \in \BZ^{|Q_0|}$ such that $\chi(d_v,e) = e(v)$ for all $e \in \BZ^{|Q_0|}$, in which case: \[
    H_*(\overline{\CM};\BC)\supset H_*(\CM;\BC) = \xi^{-1} \left ( \bigcap_{v \in Q_0} C(d_v)\right )
    \]
\end{corollary}
\begin{proof}
    The inclusion $\iota: \CM \hookrightarrow \overline{\CM}$ induces by pullback a surjection on cohomology $\iota^* : H^*(\overline{\CM};\BC) \twoheadrightarrow H^*(\CM;\BC)$ with kernel on $\CM_d$ for $d \in \BN^{|Q_0|}$ generated by relations $c_k(\mathscr{V}_v) =0$ for $k > d(v)$ (see 5.7-5.9 in \cite{Joyce2019}). Dually, $\iota$ induces an inclusion: \[
    \iota_* : H_*(\CM;\BC) \hookrightarrow H_*(\overline{\CM};\BC)
    \] and a homology class $\eta \in H_*(\overline{\CM}_d;\BC)$ lies in $H_*({\CM};\BC)$ if and only if $ c_k(\mathscr{V}_v)\cap \eta = 0 $ for $k > d(v)$. This criterion works for all $d \in \BZ^{|Q_0|}$ (dimension vectors with a negative entry will vacuously fail). In particular, if one can find $d_v$'s with $\chi(d_v,e)=e(v)$ for all $e \in \BZ^{|Q_0|}$, then the result follows. \\

    The condition that $Q$ has no oriented cycles is enough to guarantee their existence: Under a topological ordering on vertices, we see that the adjacency matrix $A$ of our quiver is strictly upper triangular, implying that $A$ is nilpotent. In particular, we can find $N > 0$ such that: \[
    (I-A)^{-1} = \sum_{k=0}^N A^k
    \] and hence $(I-A)^{-1}_{vw}$ is the well-defined number of paths from $v$ to $w$ (counting the vacuous path if $v=w$). Define $d_v$ as the $v$'th row of $(I-A)^{-1}$, so that $d_v(w)$ is the number of paths from $v$ to $w$ counting the vacuous path. It then follows that $d_v^T \cdot (I-A) \cdot e = e(v)$, as desired. 
\end{proof}
\begin{remark}
    The statement might fail if $Q$ fails to be acyclic. For example, if $Q$ only consists of a disjoint collection of vertices with no edges between them, and a self edge at each vertex, then $\chi \equiv 0$ so $C(\alpha) = V_{L_Q}$ for all dimension vectors $\alpha \in \BZ^{|Q_0|}$. 
\end{remark}

\begin{example}
    Suppose $Q$ is the $A_1$-quiver, i.e $Q_0 = \{v\}$ and $Q_1 =\emptyset$. The lattice $L_Q$ is then of rank $1$ given by $L_Q = \BZ \cdot \alpha$ where $d_v= \alpha = \textbf{1}_v$, and pairing given by $(\alpha,\alpha)=\chi(\textbf{1}_v,\textbf{1}_v)=1$. \Cref{cor: moduli-stack-of-reps-is-commutant} then tells us that the homology of the moduli stack of $\BC$-vector spaces lives inside $V_{L_Q}$ as the commutant $C(\alpha)$. 

    If $W = \langle \ket \alpha \rangle$ is the field subalgebra generated by $\ket \alpha$, then one sees that $W \subseteq C(\alpha)$ as $(\alpha,\alpha)=1 > 0$. This inclusion is actually an equality: One way to see this is to use boson-fermion correspondence, see \cite{FrenkelBenZvi2004}. Alternatively, $V_{L_Q}$ admits a bigrading by charge $n$ on $\pi_{n\alpha}$ and by weight $j_1+\cdots+j_k+\frac{n^2}{2}$ for $\alpha_{-j_1}\cdots \alpha_{-j_k} \ket{n\alpha}$, and $\ket \alpha_{(m)}$ respects this bigrading. Hence it's enough to compare bigraded dimensions. 

    Now $\xi$ identifies $C(\alpha) \cap \pi_{n\alpha}$ with $H_*(B\GL_n(\BC);\BC)$ sending homological degree $2k$ to weight $k+\frac{n^2}{2}$. Since $H^*(B\GL_n(\BC);\BZ)$ is generated by Chern classes $c_i$ of degree $2i$, the Poincaré series of $B\GL_n$ is given by: \[\sum_{j\geq 0} \dim(H^{2j}(BGL_n))q^j=\prod_{i=1}^n \frac{1}{1-q^i}=:1/(q)_n\]
    This agrees with the bigraded dimension of $W$ given by $\sum_{n \ge 0} x^n \frac{q^{n^2/2}}{(q)_n}$, see Theorem 2.6.1 in \cite{Penn2011}. 
\end{example}

\subsection{Further directions} We conclude this section by discussing questions that arose throughout this project and potential further directions. 

The usual approach to construct a lattice vertex algebra from a non-symmetric pairing $\chi$ comes by a symmetrization of said pairing $\chi^\text{sym} := \chi + \Yop{\chi}$. 
\begin{question}
    Is there a natural ``symmetrization'' functor from field algebras to vertex algebras sending the field algebra associated to $(L,\chi)$ with the vertex algebra associated to $(L,\chi^\text{sym})$?
\end{question}

One would want this functor to recover (up to isomorphism) the symmetrization for lattice field/vertex algebras.

\begin{question}
    Is there a geometric meaning to the (strong) commutant?
\end{question}

\bibliographystyle{amsalpha}               
\bibliography{refs}

\end{document}